\documentclass[a4paper,12pt,reqno]{amsart}

\usepackage{graphics}

\usepackage{soul}

\def\Infi{{\rm Inf}}
\usepackage{DLdefRS}

\graphicspath{{./QL_90x25/}}

\begin{document}

 \title[THE POINCAR\'E SERIES OF THE HYPERBOLIC COXETER
 GROUPS]{THE POINCAR\'E SERIES OF THE HYPERBOLIC COXETER
 GROUPS WITH FINITE VOLUME OF FUNDAMENTAL DOMAINS}

\author{Maxim Chapovalov, Dimitry Leites, Rafael Stekolshchik}

\date{}

\begin{abstract}
The discrete group generated by reflections of the sphere, or
Euclidean space, or hyperbolic space are said to be Coxeter groups
of, respectively, spherical, or Euclidean, or hyperbolic type. The
hyperbolic Coxeter groups are said to be (quasi-)Lann\'er if the
tiles covering the space are of finite volume and all (resp. some of
them) are compact. For any Coxeter group stratified by the length of
its elements, the Poincar\'e series (a.k.a. growth function) is the
generating function of the cardinalities of sets of elements of
equal length. Solomon established that, for ANY Coxeter group, its
Poincar\'e series is a rational function with zeros somewhere on the
unit circle centered at the origin, and gave a recurrence formula.
The explicit expression of the Poincar\'e series was known for the
spherical and Euclidean Coxeter groups, and 3-generated Coxeter
groups, and (with mistakes) Lann\'er groups. Here we give a lucid
description of the numerator of the Poincar\'e series of any Coxeter
group, and denominators for each (quasi-)Lann\'er group, and review
the scene. We give an interpretation of some coefficients of the
denominator of the Poincar\'e series. The non-real poles behave as
in Enestr\"om's theorem (lie in a narrow annulus) though the
coefficients of the denominators do not satisfy theorem's
requirements.
\end{abstract}

\maketitle

\section{Introduction}\texttt{}

The Coxeter groups split into the three types: spherical, Euclidean,
and hyperbolic. These groups are discrete reflection groups acting
on, respectively, the sphere, Euclidean space, and Lobachevsky (or
hyperbolic) space. If a hyperbolic group divides the space into
simplexes of finite volume, it is said to be of {\it Lann\'er} type
if it acts cocompactly, and of {\it quasi-Lann\'er} type otherwise.
Vinberg suggested the term in honor of Lann\'er \cite{La} who was
the first, it seems (see also \cite{CW}), to list all connected
Lann\'er diagrams (i.e., Coxeter diagrams of Lann\'er type groups);
Shwartsman and Vinberg \cite{VSh} listed all quasi-Lann\'er
diagrams.

Except for the spherical Coxeter groups $I_2^{(m)}$ (for $m\neq
3,4,6$), $H_3$, and $H_4$, each spherical (resp. Euclidean) Coxeter
group serves as the Weyl group $W_{\fg(A)}$ of simple finite
dimensional (resp. affine Kac-Moody) Lie algebra. The hyperbolic
groups of (quasi-)Lann\'er type serve as the Weyl groups of what we
suggest to call {\em almost affine} Lie algebra\footnote{These Lie
algebras are currently known under several lame names: \lq\lq
hyperbolic" (also applied to Lorentzian Lie algebras which
constitute a different set) as well as under a misleading name {\it
overextended} (it is the Dynkin diagrams that are extended twice,
not the Lie algebras). The adjective \lq\lq hyperbolic" meaningful
in the case of Coxeter groups (and helpful, unless we remember that
ALL subgroups of $O(p,1)$ are hyperbolic, while we are speaking now
only about discrete ones) is ill advised in the case of {\it these}
Lie algebras.} $\fg(A)$, where $A$ is a Cartan matrix; for the list
of almost affine Lie algebras, see the \texttt{arXiv} version of
\cite{CCLL}. We assume that all Cartan and Coxeter matrices are
indecomposable, unless otherwise stated.

\subsection{The three known facts and related problems.}
The Poincar\'e series
of the Coxeter groups of spherical and Euclidean types were known.
In this paper we explicitly compute the Poincar\'e series of certain
particular Coxeter groups of hyperbolic\footnote{The groups of
spherical and Euclidean types are often said to be of
\emph{elliptic} and \emph{parabolic} types, respectively, see
\cite{Bou,VSh}.} types.
\begin{equation}\label{F1}
\begin{minipage}[c]{14cm}
{\bf Fact 1}. {\sl Among the Coxeter groups $W$, the eigenvalues of
the Coxeter transformation of $W$ lie on the unit circle $C$
centered at the origin only for spherical or Euclidean groups
(\cite{St}). For the Coxeter groups of spherical and Euclidean
types, the zeros of the Poincar\'e series $W_{\fg(A)}$ are described
in terms of the above mentioned eigenvalues, or rather their
exponents, see Table \ref{table_grfunell}.}
\end{minipage}
\end{equation}

Our results show that for the (quasi-)Lann\'er groups (and, most
probably for all hyperbolic Coxeter groups), the zeros of the
Poincar\'e series (which, as we will show, are easy to compute) have
nothing to do with the eigenvalues of the Coxeter transformation
(which, moreover, are not easy to describe in these cases, see
\cite{St}).

\begin{equation}\label{F2}
\begin{minipage}[l]{14cm}
{\bf Fact 2}. {\sl The Poincar\'e series $W(t)$ is a rational
function for ANY infinite Coxeter group $(W, S)$ with finite set of
generators $S$. The zeros of $W(t)$ lie on on the unit circle $C$
centered at the origin, but {\bf their precise values are known only
in the spherical and Euclidean cases}. How to determine the precise
values of zeros in the other cases was unknown. The Poincar\'e
series of the Coxeter groups of hyperbolic type is exponential, so
there is a pole outside $C$ and {\bf this is all that was known
about the poles in general}.}
\end{minipage}
\end{equation}
In \cite{So,Ste,Bou}, a {\bf somewhat implicit} recurrence
expression \eqref{exBu} for $W(t)$ is given. From \cite{So, Ste,
Bou} nothing is clear about the \so{zeros of the denominators}. For
the Coxeter groups of other than spherical and Euclidean types, the
eigenvalues of the Coxeter transformations do not lie on $C$. We
show that, nevertheless,
\begin{equation*}
\begin{minipage}[l]{14cm}
\so{the zeros of the Poincar\'e series are easy to describe if these
functions are represented in a special --- virgin
--- form.}
\end{minipage}
\end{equation*}

Serre \cite{Se} was, perhaps, the first to observe several patterns
in the behavior and properties of the Poincar\'e series of the
spherical Coxeter groups:
\begin{enumerate}
  \item the Poincar\'e series are reciprocal;
  \item the value of the Poincar\'e series at $1$ is equal to the
  inverse of the Euler characteristic of the
  (geometric realization) of the respective Coxeter group.
\end{enumerate}
In works by M. Davis et. al.(\cite{DDJO, DDJMO}) the whole
Poincar\'e series, not only its value at a point, is interpreted in
terms of the weighted cohomology of Coxeter groups.

The initial goal of this note was to give an \so{explicit}
expression not only of the zeros of these rational functions (and
try to compare them with the eigenvalues of the Coxeter
transformations) but also of their {\bf poles} (not spoken about in
\cite{So,Ste,Bou} at all) in the particular cases of the {\em
(quasi-)Lann\'er groups}, i.e., Coxeter groups $(W, S)$ with
(quasi-)Lann\'er Coxeter diagrams. These groups are special in the
set of all Coxeter groups, being most close, in a sense, to the
Coxeter groups of spherical and Euclidean type: a given Coxeter
group is {\em (quasi-)Lann\'er} if its Coxeter diagram is connected,
neither spherical nor Euclidean, but any its connected proper
subdiagram is spherical (resp. spherical or Euclidean).

Knowing a recurrence formula, the problem does not seem to be
difficult ideologically but how to be sure that the result is
correct? Our own mistakes we made at first, and those we found in
the literature, make this question more serious than we thought at
first.

For the case of Coxeter diagrams with 3 vertices, see the paper by
Wagreich \cite{Wa}. Wagreich's paper is very appealing; it also
discusses several applications (e.g., due to J.~Milnor and
M.~Gromov) giving motivation for this type of activity and reasons
to publish its results in a physical journal. For applications of
Poincar\'e series of the Coxeter groups of spherical or Euclidean
type in the theory of simple finite groups, see \cite{So}. There are
other types of applications of the Poincar\'e series of the
hyperbolic groups, see, e.g., \cite{BC, GNa, DDJO}.

For the Lann\'er diagrams with 4 and 5 vertices, the answers are
known \cite{Wo}, and we used them to double check our results. We
found out that, for 5 vertices, 3 of 5 Worthington's answers are
wrong. To check our results, we need the correct results of
Worthington \cite{Wo}, and so we reproduce them. References on
Poincar\'e series of Coxeter groups include \cite{ChD, FP, Fl, Har,
Pa1, Pa2, Par}, still there is left a room to say something
reasonable.

It seemed that the denominators of the Poincar\'e series of Lann\'er
groups do not admit a nice description (and the situation with
quasi-Lann\'er groups is even worse):
\begin{equation}\label{F3}
\begin{minipage}[l]{14cm}
{\bf Fact 3}. {\sl \lq\lq With the exception of a single real pair
of poles, the poles of the Poincar\'e series of any compact
hyperbolic (Lann\'er) group with $4$ generators lie on the unit
circle $C$. This is not so for all $5$-generator Lann\'er groups"}
(\cite{CW}). 
\end{minipage}
\end{equation}
Taking the above facts into account we see the following problems:
\begin{description}\label{vopry}
  \item[(1)] Give reliable criteria for verification of the answers.
  \item[(2)] Explicitly describe the poles of the Poincar\'e series of the
5-generator Lann\'er groups.
  \item[(3)] Explicitly describe the poles of the Poincar\'e series of
quasi-Lann\'er groups.
  \item[(4)] For an infinite Coxeter group $W$, let $e(W)$, called the {\it
growth exponent}, be the inverse of the radius of convergence $R(W)$
of the Poincar\'e series $W(t)$. Compute $e(W)$, cf. \cite{Fl}.
\end{description}

\subsubsection{Our results.} We give an explicit form of
the Poincar\'e series (a.k.a. growth functions) of the Lann\'er
groups (see Tables \ref{lanner_41_44} -- \ref{lanner_51}) and
quasi-Lann\'er groups (see Tables \ref{qlanner_ql41_ql49} --
\ref{q_lanner_9_10a}).

We offer reliable means for verifications of the correctness of the
Poincar\'e series we list.

We give an interpretation of the highest and the second highest
coefficients of the denominator of the Poincar\'e series in terms of
infinite special Coxeter subgroups $W_X$ of $W$. Let $D$ be the set
of vertices corresponding to the Coxeter diagram of of the group
$W$, and $X$ run over all subsets of $D$ such that the special
subgroup $W_X$ is infinite. We say that a subgroup $W_J$ of the
Coxeter group $(W, S)$ is {\it special} (cf. \cite{Br}, p. 26) if it
is generated by a subset $J \subset S$.\footnote{In some works such
a group is called {\it parabolic}, but in other works the {\it
parabolic group} means that $w{W_J}w^{-1}$ for some $w \in W$, where
$W_J$ is the subgroup generated by $J \subset S$. Besides, the term
{\it parabolic group} is already occupied in the Lie group theory.
On top of this, instead of saying Coxeter groups of spherical and
Euclidean type some say elliptic and {\it parabolic} type,
respectively, so the term is overused, although in this context it
rhymes with hyperbolic.} The following statements hold:

 \begin{Theorem}[Theorem on the highest coefficient (Theorem \ref{th_Highest},
 $1)$] For any infinite special Coxeter group,
the coefficient $b_n$ of the highest term of the denominator of the
Poincar\'e series  $W(t)$ is as follows:
\begin{equation*}
   b_n = (-1)^{|D| + 1} - \sum\limits_{|W_X| =
   \infty}(-1)^{|X|}.
 \end{equation*}
 \end{Theorem}

\begin{Theorem}[Theorem on the second highest coefficient (Theorem
\ref{th_b_n-1})] Let $m$ be the number of factors $[n_i]$ in
  the $[n]$-complete form of the numerator of the
the Poincar\'e series  $W(t)$, see \S\ref{n-comp}. Then the
coefficient $b_{n-1}$ of the
 second highest term of the denominator is as follows:
\begin{equation*}
 b_{n-1} - m b_n  = (-1)^{|D|}|D| +  \sum\limits_{|W_X| =
 \infty}(-1)^{|X|}|X|.
 \end{equation*}
\end{Theorem}

 We derive from these theorems the following

\begin{Corollary} [On the Coxeter group with a single infinite
special subgroup (Corollary \ref{cor_sindle_inf_gr})] $1)$ For any
Coxeter group with a single infinite special subgroup, we have:
\begin{equation*}
 b_n = 0, \quad \text{ and } \quad b_{n-1} \neq 0.
\end{equation*}
In this case,  $\deg R - \deg S=1$.

$2)$ For any quasi-Lann\'er group (and also for any $1$-terminal
Coxeter group), the difference of degrees of the numerator and
denominator of the Poincar\'e series is $\deg R - \deg S\leq 1$.
\end{Corollary}

For the Euler characteristics, we have the following statement:

\begin{Proposition} [On the Euler characteristics
 (Proposition \ref{prop_Euler})]
The Euler characteristics $\chi(W)$ of the group $W$ vanishes
(equivalently, the denominator of the Poincar\'e series has the root
$t = 1$) in the following cases:

$1)$ For any affine Coxeter group.

$2)$ For any infinite (non-affine) Coxeter group $(W, S)$ with $|S|$
even. (Of course, this case includes (quasi-)Lann\'er groups whose
Coxeter diagrams have even number of vertices.)
\end{Proposition}

We have found out that the poles of the Poincar\'e series of the
quasi-Lann\'er groups behave rather nicely:

\ssec{Towards a generalization of the Enestr\"om theorem}
\sssec{Gal's formulation} For recent studies of the poles of the
Poincar\'e series of Coxeter groups, see Gal's interesting preprint
\cite{Gal} with preliminary results of an aborted research. Gal
considered Coxeter diagrams for which the nerve $N_W$ (see subsec.
\ref{Nerve}) of the corresponding Coxeter group $W$ is a homology
sphere\footnote{A {\it homology sphere} is an $n$-dimensional
manifold having the same homology groups as $S^n$ does.}.

Gal wondered how many real poles can the Poincar\'e series of such a
group have (he notes that the degree of the denominator of the
Poincar\'e series of any non-right-angled Coxeter group may be
however greater than the dimension of its nerve). If $W$ is an
affine Coxeter group, then there is a unique real pole of order $n$
at $1$ \cite{Bou}. If $\dim N_W=n \leq 3$, then there are exactly
$n$ positive real poles \cite{Par}. Moreover, in these two cases,
all the non-real poles lie on the unit circle.

Gal writes that usually (but does not explain what the share of this
\lq\lq usually" in the general picture is and what the exceptions
are), if $\dim N_W\geq 3$, the non-real poles of the Poincar\'e
series fail to lie on the unit circle. Looking at the examples known
to him Gal made the following observation (he writes that he \lq\lq
tested a number of groups whose nerve is a simplex or a product of
simplexes" but, regrettably, did not specify the number and gave
only two illustrations which, actually, are $L5_5$ and $QL10_2$):
\begin{equation}\label{poles}
\text{\begin{minipage}[c]{10cm} {\sl several poles lie \lq\lq
near\rq\rq  the real positive half-line and the rest of the poles
tend to lie \lq\lq near\rq\rq  the unit circle.}\end{minipage}}
\end{equation}

We do not know how to quickly say if the nerve of $W$ is a homologic
sphere or not, but the examples Gal gives made us wonder if not just
two but ALL the cases we study satisfy \eqref{poles}. Indeed, they
are, with several corrections of Gal's description:

\sssec{Quasi-Lann\'er case. General hyperbolic Coxeter groups}
Having found the precise expressions of the Poincar\'e series and
their poles we saw that the distribution of poles, which could have
been random, does resemble the pattern \eqref{poles} almost
correctly described by Gal \cite{Gal}. Let us forget for a moment
the poles lying \lq\lq near the real positive half-line"; the
remaining poles do lie in a thin annulus concentric with and
sometimes containing the unit circle.

Our results and Gal's hints lead us to a result of G.~Enestr\"om
\cite{E}. His theorem (rediscovered by Kakeya \cite{Kak}, see
interesting reviews \cite{GH, vV} and references therein; Kakeya's
work had some mistakes but, despite this, the statement is often
referred to as Enestr\"om-Kakeya theorem) says

\sssbegin{Theorem}\label{EK} Let $p(t)=a_0+a_1t+\dots +a_nt^n$ be a
polynomial with positive coefficients, set
$m:=\mathop{\min}\limits_{0\leq i< n}\frac{a_i}{a_{i+1}}$, and
$M:=\mathop{\max}\limits_{0\leq i< n}\frac{a_i}{a_{i+1}}$. Then all
the roots of $p(t)$ lie in an annulus with bounding circles of
radius $m$ and $M$ concentric with the unit circle $C$ centered at
the origin.
\end{Theorem}
The coefficients of the denominators of the Poincar\'e series of the
(quasi-)Lann\'er groups do not satisfy the conditions of the
Enestr\"om theorem but the \so{non-real} zeros of these polynomials
behave as if they do, or almost: all non-real roots lie in an
annulus centered at the origin (except that we do not know how to
define the radii $m$ and $M$ of the bounding circle from the
coefficients and the annulus does not necessarily contain $C$).

It is natural, therefore, to disregard for a moment the real roots
and try to find the conditions the coefficients of the denominators
of the Poincar\'e series of the (quasi-)Lann\'er groups satisfy in
order to derive a generalization of the Enestr\"om theorem for
polynomials whose real coefficients can be of any sign or vanish.

At our request, V.~Molotkov studied several simplest Lann\'er cases
and saw that the poles lying on $C$ are hardly roots of unity
(unlike the zeros of the numerators of the Poincar\'e series of all
Coxeter groups). He also observed that, in contradistinction with
what is depicted in Gal's illustration for $QL10_2$, none of the
non-real poles is lying \lq\lq near\rq\rq the real positive
half-line \lq\lq parallel to it\rq\rq. Instead
\begin{equation}\label{poles1}
\text{\begin{minipage}[c]{13cm} \so{the non-real roots lie in a {\bf
thin} annulus concentric with the unit circle $C$; the real poles
(if any) lie near $1$ or $-1$}.
\end{minipage}}
\end{equation}
Molotkov's results, more precise than Gal's, inspired us to verify
and sharpen Gal's conjecture \eqref{poles} as formulated in
\eqref{poles1}; in most cases, NONE of the non-real roots lies on
$C$. Bar few exceptions for 4-vertex diagrams, the poles we found
numerically are non-simple-looking (for humans) algebraic numbers.
Therefore we have summarized the answer by listing only the real
roots and the extremal values of the absolute values of the non-real
roots, see Tables \ref{Tabpoles}--\ref{Tabpoles10}. We conjectured
that \so{the non-real poles of the Poincar\'e series of any Coxeter
group $(W, S)$ with $|S|<\infty$ lie in a thin annulus}: This was
the case with several of the Coxeter groups we unwillingly
considered while making typos in the input data. However, we tested
the conjecture on the {\bf reflective arithmetic} Coxeter groups
(\cite{VSh}, Table in subsect. 2.1) and {\bf non-arithmetic} Coxeter
groups (\cite{VSh}, Table in subsect. 3.2) and found out that this
conjecture is overoptimistic: \so{Most} the non-real poles of the
Poincar\'e series of these Coxeter groups lie in a thin annulus, but
not all.

\sssbegin{Problem} What are the conditions on the coefficients of
the real polynomial for its non-real roots to lie in a thin annulus?
How to describe the radii of the circles that bound the annulus in
terms of the coefficients of the polynomial?
\end{Problem}

\ssec{Discussion: Infinite Coxeter groups} We conclude from the
results of the paper that amount and interrelation of infinite
special subgroups in the given infinite Coxeter group is very
essential and closely related to predicting coefficients of the
Poincar\'e series. This motivated us to divide the set of all
infinite Coxeter groups as follows.

\sssec{$k$-terminal Coxeter groups}

 We say that an infinite Coxeter group $W$ is {\bf
$k$-terminal} if the length of any chain of its infinite special
subgroups ordered by inclusion is $\leq k$ and at least one chain is
of length $k$. Examples:

\begin{equation}
\label{k_term}
\begin{split}
\{\text{0-terminal} & \text{ Coxeter groups} \} = \{ \text{affine
Coxeter groups} \}
   \quad \cup \quad
   \{\text{Lann\'er Coxeter groups} \} \\
     \{\text{1-terminal} & \text{ Coxeter groups}\}  = \{\text{quasi-Lann\'er Coxeter groups}
   \}
\quad \cup \quad \\
  & \{\text{Coxeter groups with all special subgroups finite, affine or Lann\'er} \}
   \end{split}
\end{equation}

Let $W$ be any $k$-terminal Coxeter group and $\mathcal{P}(W)$ the
poset of all infinite special subgroups of $W$. Let $l(E)$ be the
{\it level of an element} $E \in \mathcal{P}(W)$ defined so that
$l(W) = 0$, $l(E) = 1$ for any maximal infinite special subgroup,
and so on. Denote by $\Infi_m$ the number of infinite special
subgroups of level $m$ in the poset $\mathcal{P}(W)$.

\sssec{Subsets $I^n_i$ of infinite Coxeter groups} Set

$I^{F}$ :=  \{ finite Coxeter groups \},

$I^{0}_1$ := \{ affine Coxeter groups \},

$I^{0}_2$ :=  \{ Lann\'er Coxeter groups \}, \quad $I^{0} = I^{0}_1
\cup I^{0}_2$,

$I^{1}_1$ := \{ quasi-Lann\'er Coxeter groups \} := \{ $W \mid$
every proper special subgroup of $W' \subset W$ is a group from
$I^{F} \cup I^{0}_1$, and there exists $W' \subset W$ such that
 $W' \in I^{0}_1$ \},

$I^{1}_2$ :=  \{ $W \mid$ every proper special subgroup
 $W' \subset W$ is a group from $I^{F} \cup I^{0}$, and there exists
 $W' \subset W$ such that  $W' \in I^{0}_2$\}, \quad $I^{1} = I^{1}_1  \cup
I^{1}_2$,

Let us introduce by induction subsets $I^{n}_1$, $I^{n}_2$ and
$I^{n}$ as follows:

$I^{n}_1$ := \{ $W \mid$ every proper special subgroup
 $W' \subset W$ is a group from $I^{F} \cup ( \bigcup\limits_{i = 0}^{n-2}I^{i})
\cup I^{n-1}_1$, and there exists $W' \subset W$ such that
 $W' \in I^{n-1}_1$  \},

$I^{n}_2$ :=  \{ $W \mid$ every proper special subgroup
 $W' \subset W$ is a group from $I^{F} \cup (\bigcup\limits_{i = 0}^{n-1}I^{i})$,
and there exists $W' \subset W$ such that  $W' \in I^{n-1}_2$ \},
\quad $I^{n} = I^{n}_1 \cup I^{n}_2$.

\sssbegin{Proposition} For $n \geq 0$, we have

$1)$ the set $I^{n}$ consists of $n$-terminal Coxeter groups,

$2)$ $I^{n}_1 \cap I^{n}_2 = \emptyset$.
\end{Proposition}

The proposition is easy to prove by induction. \qed.

\sssec{Classification problem}  As the next natural step in the
study of infinite Coxeter groups,  it seems to us important to
describe the set $I^{1}_2$, the next after quasi-Lann\'er Coxeter
groups in hierarchy of
 the $k$-terminal Coxeter groups.

\section{Precise setting of the problems}
\ssec{Generating functions} Generating functions of graded objects
were introduced and studied by Hilbert and Poincar\'e at more or
less the same time. Leaving touchy priority questions aside,
Wikipedia informs us:
\begin{equation}\label{HPser}
\text{\begin{minipage}[c]{14cm} \lq\lq A {\it Hilbert-Poincar\'e
series}, named after David Hilbert and Henri Poincar\'e, is an
adaptation of the notion of dimension to the context of graded
algebraic structures (where the dimension of the entire structure is
often infinite). It is a formal power series in one indeterminate,
say $t$, where the coefficient of $t^n$ gives the dimension (or
rank) of the sub-structure of elements homogeneous of degree $n$."
\end{minipage}}
\end{equation}

Observe that in the above definition certain restrictions are taken
for granted: the dimension of each homogeneous component must be
finite, and only non-negative components are usually non-zero;
\lq\lq graded" is only assumed to be by means of $\Zee$; for
$\Zee^k$-graded objects (under similar restrictions: The support of
the degrees with non-empty components lies in the cone with
non-negative coordinates and each component is finite), we get
series in several indeterminates, as in \cite{McD, DDJO} and Table
\ref{table_grfunell}.

In the particular case of Coxeter groups stratified by the length of
their elements, instead of the term \lq\lq Hilbert-Poincar\'e
series" the term {\it Poincar\'e series} is usually used, and lately
it is called also by the {\it growth function}. We use the
terminology of the classics, and in a particular case of Coxeter
groups of (quasi-)Lann\'er type is the object of our study.

\subsection{Coxeter groups.}
 A {\it Coxeter system} to be a pair $(W, S)$ consisting of a group
 $W$ and a set of generators $S \subset W$ subject to relations
\begin{equation}
 \label{coxgr}
 (st)^{m_{s, t}}=1,\text{ where }
   m_{s, s} = 1, m_{s, t}=m_{t, s}\geq 2  \text{ for } s \neq t \text{ in } S.
\end{equation}
If no relation occurs for a pair $s, t$, then it is assumed that
$m_{s, t}=\infty$. In this presentation $W$ is a {\it Coxeter
group}. The symmetric matrix $M=(m_{s, t})_{s, t\in S}$ is called a
{\it Coxeter matrix}.

The presentation of every finitely generated Coxeter group can be
illustrated by an undirected labeled graph, called {\it Coxeter
diagram}, whose vertices correspond to the generators $S$ of $W$ and
edges are as follows. If $m_{s, t}=2$ then no edge joins $s$ and
$t$. If $m_{s, t}=3$, then an edge joins $s$ and $t$. The edge
between the vertices corresponding to $s, t\in S$ is endowed with
label $m_{s, t}$ if $m_{s, t}>3$.

The Poincar\'e series $W_{W, S}(t)$ of a group $W$ relative to a
finite generating set $S$ is briefly denoted $W(t)$
 and defined as follows. For any $g\in W$, define the {\it length}
$l(g)$ to be the minimum length of all words in $S$ representing
$g\neq 1$ and $l(1)=0$. Then
\\
\vspace{-1mm}
 \begin{equation}
  \label{GSt}
  W(t) := \sum_{g\in W} t^{l(g)}.
 \end{equation}

\subsubsection{Remarks}
1) The Coxeter diagrams, so graphic for Weyl groups of finite
dimensional and Kac-Moody Lie algebras, are utterly useless if the
Coxeter matrix is not sparse, as is the case of Lorentzian Lie
algebras considered by Borcherds, and Gritsenko and Nikulin, see
\cite{GN}, or in the cases considered in sect. \ref{Floyd}. In this
note, we deal with the cases where graphs are helpful, but the
reader should realize that {\it actually} we deal with Coxeter
matrices.

2) Other notation used (less convenient, we think, if there are many
cases of multiple edges): The edge between nodes $s$ and $t$ is
often depicted as a multiple one of multiplicity $m_{s, t}-2$,
unless $m_{s, t}=\infty$; for $m_{s, t}=\infty$, the edge is usually
depicted {\bf thick}.

For the Lie algebra $\fg(A)$ with Cartan matrix $A$ normalized, as
usual, so that $A_{ii}=2$, and with non-positive integer
off-diagonal elements, the Coxeter matrix $M=(m_{ij})_{ij\in S}$ is
given by the conditions
\begin{equation}
 \label{cos} \begin{tabular}{||c|c|c|c|c|c||}\hline\hline
 $A_{ij} A_{ji}$&0&1&2&3&$\geq 4$\cr
 \hline
 $m_{ij}$&2&3&4&6&$\infty$\cr\hline\hline
\end{tabular}
\end{equation}
\\
We do not reproduce the list of spherical (resp. Euclidean) Coxeter
diagrams (see \cite{Vi}): They are easily obtained from the
well-known Dynkin graphs and their Cartan matrices, see \cite{Bou},
(resp. from their extended versions, see \cite{K, St}).

\ssec{Exponents} Let $W$ be a finite group generated by reflections
$r_i$, where $i=1,\ldots, n$, in the Euclidean space or,
equivalently, on the sphere. (For example, the Weyl group
 $W = W_{\fg}$ of a simple Lie algebra $\fg$ naturally acts in the root
space of $\fg$.) Let ${\bf C}:=\prod r_i$ be the product of all
generators (in any order; all these products are conjugate, see
\cite{St}). For the Weyl groups of simple finite dimensional and
affine Kac-Moody Lie algebras, the eigenvalues of ${\bf C}$ are of
the form $\omega^{m_i}$, where $\omega=e^{2\pi i/h}$ and where
$h=1+\max m_i$ is the Coxeter number
--- the order of ${\bf C}$ (\cite{CM}, \cite{OV}, \cite{St}). The
numbers $m_i$ are called the {\it exponents} of the corresponding
Coxeter group, see \cite[Table 2]{Cox}, and our Table
\ref{table_expon}.

Here is an excerpt from \cite[p.765]{Cox} regarding exponents (at
places in our own words):
\begin{equation}
 \label{expo}
\begin{minipage}[l]{14cm}
\lq\lq Most of the applications of ${\bf C}$ are related with $h$.
We consider the characteristic roots
\[\omega^{m_1}, \ldots, \omega^{m_n}\]
of ${\bf C}$ and the exponents are certain integers which may be
taken to lie between $0$ and $h$. They are computed by a
trigonometrical formula involving the periods [i.e., orders]
$m_{ij}$ of the products of pairs of generators. (The product of two
reflections is simply a rotation.)

The point of interest is that the same integers occur in a different
connection. It turns out that the order of the group is
\begin{equation*}
 (m_l+ 1)(m_2 + 1)\ldots (m_n + 1),
\end{equation*} and that these factors $m_i + 1$ are the degrees of $n$ basic
invariant forms [William Burnside, {\em Theory of Groups of Finite
Order}, Cambridge, 1911; Chapter XVII]. Moreover, when every
$m_{ij}$ is equal to 2, 3, 4 or 6, so that the group is
crystallographic, there is a corresponding continuous group $G$, and
the Betti numbers of the group manifold are the coefficients in the
Poincar\'e polynomial (of the manifold of the Lie group $G$)
\[(1+t^{2m_l+ 1})(1+t^{2m_2 + 1})\ldots (1+t^{2m_n + 1})."\]
\end{minipage}
\end{equation}

\subsection{The Poincar\'e series (a.k.a. growth functions) of the
Coxeter groups.} Following Solomon, Bourbaki \cite{Bou} gives an
explicit expression of the Poincar\'e series $W_{\fg}(t)$ for the
Weyl group $W_{\fg}$ of simple finite dimensional Lie algebra $\fg$
in terms of exponents:
\begin{equation}
 \label{solo}
W_{\fg}(t) = \prod\frac{1-t^{m_i+1}}{1-t}.
\end{equation}
This formula is applicable not only to the Weyl groups of the simple
finite dimensional Lie algebras, but to other groups of spherical
 type, see Table \ref{table_grfunell}.

The generalization of \eqref{solo} to affine Weyl groups is due to
Bott \cite{Bo}; see also Reiner's notes \cite{Rei} with an
exposition of the proof of Bott's result due to Steinberg
\cite{Ste}. Bott keeps writing about the loop groups or loop
algebras (i.e., algebras of the form $\tilde\fg:=\fg\otimes
\Cee[u^{-1}, u]$, where $\fg$ is any simple finite dimensional Lie
algebra), but in reality he only considers the Weyl groups of the
Lie algebras of these loop groups. Since the exponents are defined
up to dualization of the root system, the Poincar\'e series for the
\lq\lq twisted" affine Kac-Moody algebras are covered by Bott's
result. The answer is given by the formula
\begin{equation}
 \label{HPWloopquot}
W_{\tilde\fg}(t) = \prod\frac{1-t^{m_i+1}}{(1-t)(1-t^{m_i})} =
 W_{\fg}(t)\prod\frac{1}{1-t^{m_i}}.
\end{equation}
Let us now try to perform the next step
--- consider the Weyl groups of almost affine Lie algebras.

\footnotesize

\subsection{Digression: (Quasi-)Lann\'er groups are the Weyl groups of
almost affine Lie algebras.} There are several (intersecting but
distinct) sets of Lie algebras whose elements are often called
\lq\lq hyperbolic" Lie algebras. We would like to carefully
distinguish between these sets so need an appropriate name for each.
We say that a submatrix of a square matrix is {\it principal} if it
is obtained by striking out a row and column that intersect on the
main diagonal. We say that Lie algebra with Cartan matrix whose
entries belong to the ground field is {\it almost affine} if it is
not finite dimensional or affine, and its subalgebra corresponding
to any principal submatrix of the Cartan matrix is the sum of finite
dimensional or affine Lie algebras.

Z.~Kobayashi and J.~Morita classified the almost affine Lie algebras
with indecomposable symmetrizable Cartan matrix of size $>2$
\cite{KoMo}. Later, Li Wang Lai \cite{Li} obtained a complete answer
(for Cartan matrices of size $>2$): there are 238 almost affine Lie
algebras; 142 of these algebras have a symmetrizable Cartan matrix.
Later Sa{\c{c}}lio{\u{g}}lu \cite{S} rediscovered the result of
Kobayashi and Morita (with few omissions); his paper is devoted to
physical applications and is very interesting.

In this paper we derive explicit formulas for the Poincar\'e series
of the groups most close in a sense to the Weyl groups of simple
finite dimensional Lie algebras. In the literature, in similar
studies, the authors write sometimes that they are studying the Lie
algebras or even the Lie groups having these Lie algebras, whereas
they are only studying the Weyl groups of these Lie algebras. This
subtlety is sometimes important: In particular, to list all the
groups we are dealing with (Lann\'er and quasi-Lann\'er) is much
easier than to list the Lie algebras whose Weyl groups they are.
These are almost affine (a.k.a hyperbolic) Lie algebras; their
complete list was unknown when the description of the growth
functions of their Weyl groups has begun (and the classification of
these Lie algebras is not needed in this particular study of their
Weyl groups). There are several stages of generalization of simple
finite dimensional Lie algebras (which all possess very particular
Cartan matrices) to the Lie algebras with more-or-less arbitrary
Cartan matrix. We intend to generalize the results on the growth
functions known for the Weyl groups of simple finite dimensional and
affine Kac-Moody Lie algebras to the case of Weyl groups of almost
affine Lie algebras. These Lie algebras became of acute interest
lately in connection with \lq\lq cosmic billiards"; for details and
further references, see \cite{H}, \cite{BS}. The Poincar\'e series
of the Weyl groups of almost affine Lie algebras are invariants of
these Lie algebras that can be used further, see \cite{Wa} and
references therein. The set of almost affine Lie algebras has a
non-empty intersection with the (different) set of Lorentzian Lie
algebras, sometimes also called \lq\lq hyperbolic". For applications
of Lorentzian Lie algebras, see \cite{RU}, \cite{GN}. For one of
these applications Borcherds was awarded with Fields medal.

\normalsize

\section{The Poincar\'e series (known facts)}
\subsection{The Solomon-Steinberg recursion \eqref{exBu}.}
For any finite set $X$, let $\eps(X)=(-1)^{\card(X)}$. Let $W_X(t)$
be the Poincar\'e series (a polynomial or series) of the Coxeter
group $W_X$ whose Coxeter graph is $X$. If $\card W_D<\infty$, let
$\cM$ be the maximal length of the elements of $W_D$ (there is only
one element of maximal length).

Ex. 26 to \S1 of Ch.4 \cite{Bou} claims that \so{for any Coxeter
graph $D$, we have} (this formula is obviously due to Solomon
\cite{So} (although in particular cases of finite Weyl groups this
may have been established earlier by Chevalley, see \S3.15 in
\cite{Hu}); Steinberg \cite{Ste}, Theorem 1.25 gave a simpler proof;
for an exposition of Steinberg's proof, see also \cite{McD}, where
there are considered multiparameter series \so{taking into account
difference in length of roots}\footnote{Therefore, for this task, we
need not just Coxeter graphs (i.e, Coxeter matrices) but the Dynkin
diagrams (Cartan matrices), and hence the classification of almost
affine (a.k.a. hyperbolic) Lie algebras due to \cite{Li, S}; for the
list of such diagrams/matrices, see also \cite{CCLL}.}); here $X$ is
any {\bf complete}\footnote{Recall that a  subgraph is {\it
complete} if each of its nodes is connected to every other of its
nodes.} subgraph of $D$:
\begin{equation}
 \label{exBu}
 \sum_{X\subset D} \frac{\eps(X)}{W_X(t)}=
 \begin{cases}
 \displaystyle\frac{t^{\cM}}{W_D(t)}& \text{if $\card W_D<\infty$}, \\
 0   & \text{otherwise}.
 \end{cases}
\end{equation}
In this expression, the summand corresponding to the empty subgraph
is equal to $1$.

Recall that the {\it rational} (non-polynomial) function $P(t)$  is
said to be {\it reciprocal} if $P(t^{-1}) = P(t)$; if $P(t^{-1}) =
-P(t)$ the rational function $P(t)$ is often said to be {\it
anti-reciprocal}.

The {\it polynomial} function $P(t)$ is said to be {\it reciprocal}
(resp.{\it anti-reciprocal}) if
\[P(t)=t^{\cM}P(t^{-1}), \text{~~(resp. $P(t)=-t^{\cM}P(t^{-1})$), where $\cM=\deg P$}. \]

The (anti-)reciprocal function is said to be $\pm$-{\it reciprocal}.

The recurrence \eqref{exBu} and $\pm$-reciprocity  of $W_X(t)$ if
$|W_X|<\infty$ imply the following sharpening of \eqref{exBu} due to
Steinberg \cite{Ste}: If $\card W_D=\infty$, then
\begin{equation}
 \label{exSte}
 \frac{1}{W_D(t^{-1})}=\sum_{X\subset D\;\;\mid\;\;
 \card W_X<\infty}\frac{\eps(X)}{W_X(t)}.
\end{equation}

To begin the induction, recall the following facts:

0) If the Coxeter graph $X$ is the disjoint union of its connected
components $X_i$, then $W_X(t)=\prod W_{X_i}(t)$. Hereafter it is
advisable to use standard simplified notation: For any
$n\in\Nee\cup\{\infty\}$, set
\begin{equation}
 \label{[n]}
 [n] :=
 \begin{cases}1 + t + \dots + t^{n- 1}&\text{for $n<\infty$}, \\
 1 + t + \dots = \displaystyle\frac{1}{1-t}&\text{for
$n=\infty$}.\end{cases}\end{equation}

1) $W_\emptyset(t)=1$ and $W_{*}(t)=1+t=[2]$ (that is, for the graph
consisting of 1 vertex and 0 edges).

2) If $X$ has two vertices joined by $m-2$ edges, then
\begin{equation}
 \label{PX}
 W_X(t)=\begin{cases}
 \displaystyle\frac{(1+t)(1-t^{m+2})}{1-t}=[2][m+1]&\text{if
 $3\leq m<\infty$ (for $I_2^{(m)}$)}, \\
 \displaystyle\frac{1+t}{1-t}=[2][\infty]&\text{if $m=\infty$ (for
$I_2^{(\infty)}$)}.\end{cases}
\end{equation}

3) The Poincar\'e series of the 3-generator Coxeter group $G_{p, q,
r}$ with diagram $L3$ or $QL3$ (if $|G_{p, q, r}|<\infty$, then
$\frac1p + \frac1q + \frac1r
> 1$):
\begin{equation}\label{pqrfin} W_{G_{p, q, r}}(t)=\frac{[2][p][q][r]}{
[2][p][q][r]-3[p][q][r] + [p][q] + [p][r] +
[q][r]}\times\begin{cases}(t^\cM + 1)&\text{if $|G_{p, q,
r}|<\infty$, }\\1&\text{otherwise}. \end{cases}
\end{equation} The Coxeter graphs are as
follows:

$L3$: Each diagrams on 3 vertices is a triangle with edges labeled
by $p$, $q$, $r$ such that $2\leq p, q, r< \infty$ and $\frac1p +
\frac1q + \frac1r < 1$. One (only one) of the labels $p$, $q$, $r$
may be equal to 2, and then the graph is not, actually, a triangle.

$QL3$: The graphs look as those for $L3$ but any of the labels $p$,
$q$, $r$ may be (and at least one is) equal to $\infty$.

We summarize the results needed to explicitly compute \eqref{exSte}
in Table \ref{table_grfunell}.

\ssec{Lann\'er and quasi-Lann\'er diagrams on $>3$ vertices} In the
literature we saw, these diagrams are seldom identified (the only
exception known to us is an interesting paper \cite{JKRT} with too
complicated\footnote{In addition to overcomplicated proper names,
called {\it Witt symbols}, there are given in \cite{JKRT} also {\it
Coxeter symbols} that less vaguely encode some of the Coxeter graphs
but can not be used as short names, either, and are not clearly
defined for an arbitrary diagram in either \cite{CM} or \cite{JKRT}
(try to reconstruct the rules for, e.g., $\overline{DP}_3$,
$\overline{M}_3$ or $\overline{N}_4$).} names for them), so we
simply number them for convenience. The first to list these diagrams
was, it seems, Lann\'er \cite{La}, see also \cite{CW50}.

For the Lann\'er diagrams and the corresponding Poincar\'e series,
see Tables \ref{lanner_41_44} -- \ref{lanner_51}.

For the quasi-Lann\'er diagrams and the corresponding Poincar\'e
series, see Tables \ref{qlanner_ql41_ql49} -- \ref{q_lanner_9_10a}.

\sssec{Worthington's results} For the Lann\'er groups with 4
generators, Worthington computed their Poincar\'e series, and we
confirm them in Table \ref{lanner_41_44}--\ref{lanner_46_49}.
Worthington computed the Poincar\'e series of the quasi-Lann\'er
groups with 5 generators, but in 3 of 5 cases his answers are wrong.

\ssec{Two interesting (and correct although strange) --- but useless
for us
--- formulas}\label{Floyd} Floyd and Plotnick \cite{FP} cite the
following statement they attribute to Parry. The first displayed eq.
on p.524 of \cite{FP} gives the following presentation of a Coxeter
group $W$:
\begin{equation}
 \label{PresFP}
 W=\left\langle g_1, \dots, g_d\mid g_i^2, (g_ig_{i+1})^{a_i}\right\rangle.
\end{equation}
Obviously $i$ runs 1 through $d$ and --- although this was not
mentioned ({\it sapienti sat}) --- the relation between $g_d$ and
$g_{d+1}$ should be understood as a relation between $g_d$ and
$g_{1}$. Since nothing is mentioned, there are no relations between
$g_i$ and $g_{j}$ for $i\neq j\pm 1$ (and $g_d$ and $g_{j}$ for
$j\neq d-1$ or 1), i.e., there are lots of relations of the form
$(g_ig_{j})^{\infty}=1$. Then (as usual, the hatted factor should be
ignored)
\begin{equation}
 \label{Parr}
 W(t)=\frac{[2][a_1]\dots[a_n]}{(t+1-n)[a_1]\dots[a_n]+\sum
 [a_1]\dots\widehat{[a_i]}\dots[a_n]}.
\end{equation}
In 1991, Floyd submitted a paper \cite{Fl} in which the following
condition of applicability of formula \eqref{Parr} is added to tho
conditions given in \cite{FP}: The set $a_1$, \dots, $a_{n}$ is said
to be {\it unacceptable} if $(a_1, \dots, a_{n})= (2, 2, 2, 2)$ for
$n=4$ or $\frac{1}{a_1}+\frac{1}{a_2}+\frac{1}{a_{3}}\geq 1$ for
$n=3$. For all other --- \so{acceptable} --- sets of labels, the
following formula is offered in place of \eqref{Parr}:
\begin{equation}
\label{Fl} W(t)=\frac{[2][a_1]\dots[a_n]}{[2][a_1]\dots[a_n]-t\sum
[a_1]\dots[a_i-1]\dots[a_n]}.
\end{equation}
The MAIN applicability condition of \eqref{Parr} and \eqref{Fl}
mentioned in \cite{FP, Fl} is, however, that
\begin{equation}
\label{applic} \text{the group $W$ should act on the 2-dimensional
hyperbolic space.}
\end{equation}
The mysterious (how to verify it for an abstractly given group?)
condition \eqref{applic} is applicable to the (quasi-)Lann\'er
graphs only in certain cases of three vertices, so it is of no
interest to us. We have included the remarkable
--- they are symmetric in the $a_i$ which is astounding ---
formulas \eqref{Parr} and \eqref{Fl} in this paper for completeness
of the picture.

\ssec{The multiparameter case} Macdonald \cite{McD} describes the
passage to the multiparameter case with his usual clarity:

For any Coxeter group $(W, S)$, let $S_i$ for a set of indices $I$,
be the equivalence classes of the relation \lq\lq $s$ is conjugate
to $s'$ in $W$" between elements of $s,s\in S$.
\begin{equation}
\label{McDrule}
\begin{minipage}[c]{13cm}
The subsets $S_i$ can be read off the Coxeter diagram of the group
$(W, S)$ by \lq\lq reduction modulo 2": if we delete from the
diagram all bonds bearing the even label or $\infty$, then the
connected components of the resulting graph correspond to the $S_i$.
\end{minipage}
\end{equation}

Let ${\bf t}=(t_i)_{i\in I}$ be a family of indeterminates. For each
$s\in S$, define
\begin{equation}
\label{multip} t_s:=t_i \text{~~ for $s\in S_i$}.
\end{equation}
Let $w=s_{i_1}\dots s_{i_l}$ be a reduced decomposition of any
 $w \in W$, i.e., representation of $w$ as the product of the least number
of generators. Then the monomial
\begin{equation}
\label{multip2} t_w:=t_{s_{i_1}}\dots t_{s_{i_l}}
\end{equation}
does not depend on the choice of reduced decomposition. Then,
clearly,
\begin{equation}
\label{multip3} W({\bf t}):=\sum_{w\in W}t_{w}.
\end{equation}
Let $l_i(w)$ be the $i$-length of $w$, i.e., the number of the
generators in the reduced decomposition of $w$ belonging to $S_i$.
Then
\begin{equation}
\label{multip4} t_w:=\prod_{i\in I}t_{i}^{l_i(w)}.
\end{equation}
With these definitions, the formula \eqref{exBu} is still true with
${\bf t}$ instead of $t$. For the necessary changes, see Table
\ref{table_grfunell}: For the Coxeter groups of spherical type,
$|I|\leq 2$ and $|I|= 2$ only in the three cases.

Observe several subtleties:

1) Clearly, the Coxeter diagram $D$ is not sufficient to describe
the multiparameter Poincar\'e series: We have to distinguish between
short and long roots if $|D|>2$, so we have to distinguish between
the $B_n$ and $C_n$ cases.

2) Although the Lie algebras with non-symmetrizable Cartan matrices
do have the Weyl group defined by eq. \eqref{cos}, Macdonald's rule
\eqref{McDrule} is only applicable to the root systems described by
symmetrizable Cartan matrices: Otherwise the notion of short/long
root is not well-defined.

Although Macdonald's paper is devoted to all Coxeter groups of
spherical and Euclidean cases, he evaded computing the
multiparameter Poincar\'e series for the Weyl groups of the twisted
loop Lie algebras leaving this as \lq\lq an easy task for the
reader", having indeed explained all the needed steps. This was,
perhaps, a joke: all one should do is to renumber the indeterminates
in accordance with Table \ref{table_grfunell} minding the above
subtlety 1). Macdonald missed (or left as a trivial exercise?) the
case of $A_1^{(1)}$ (the answer for which coincides with that for
$A_2^{(2)}$, see Table \ref{tab_multip}).

Unless the authors of \cite{DDJO}, where the multiparameter growth
functions are applied, or somebody else, will ask us to do the job,
we intend to imitate the behavior of Prof. Macdonald, and redirect
the reader: For the classification of Li and Sa{\c{c}}lio{\u{g}}lu,
see more accessible list in \cite{CCLL}; the code is available at
\cite{DCh} and how to proceed with the code is described in the last
section of this work. The task is now routine while to list the
results will double the length of the tables.

\section{The Euler characteristic (from [D], [Br])}

\ssec{The geometric realization of the simplicial complex}
  \label{geom_realiz}
A {\it simplicial complex} with vertex set $\cV$ is a collection
$\Delta$ of finite subsets of $\cV$ (called {\it simplexes}) such
that every singleton $\{v\}$ is a simplex and every subset of a
simplex $A$ is a simplex (called a {\it face} of $A$), \cite[Ch.I,
App.]{Br}. The cardinality $r$ of $A$ is called the {\it rank} of
$A$, and $r-1$ is called the {\it dimension} of $A$. We include the
empty set as a simplex; it has rank $0$ and dimension $-1$. A {\it
subcomplex} of $\Delta$ is a subset $\Delta'$ which contains, for
each of its elements $A$, all the faces of $A$,  thus $\Delta'$ is a
simplicial complex in its own right, with vertex set equal to some
subset of $\cV$. Note that $A$ is a poset, ordered by the face
relation.

The {\it geometric realization} $|\Delta|$ of $\Delta$ is a
topological space partitioned into (open) simplexes $A$, one for
each non-empty $A \in \Delta$. This  topological space is
constructed as follows: We start with an abstract real vector space
with $\cV$ as a basis. Let $|A|$ be the interior of the simplex in
$\cV$ spanned by the vertices of $A$, i.e., $|A|$ consists of the
linear combinations $\sum\limits_{v \in A}{\lambda_v}v$ with
$\lambda_v > 0$ for all $v\in A$ and $\sum\limits_{v \in
A}{\lambda_v} = 1$. We then set
\begin{equation*}
    |\Delta| = \bigcup\limits_{A \in \Delta} |A|.
 \end{equation*}

\ssec{Cells and chambers}

Let $V$ be the  space of a geometric {\it realization} of the
Coxeter group $(W, S)$, and $\dim V = n$. Let $\mathbb{H} =
\{H_1,... ,H_k\}$ be an arbitrary finite set of hyperplanes in $V$.
The hyperplanes $H_i$ cut $V$ into polyhedral pieces by means of
reflections $s_i:=s_{H_i}$ that generate $W$ . For each $i = 1,...,
k$, where $k = \card(S)$, let $f_i : V \tto \mathbb{R}$ be a
non-zero homogeneous linear function that singles out $H_i$ by the
equation $f_i = 0$. The function $f_i$ is uniquely determined by
$H_i$, up to a non-zero factor.

\begin{figure}[h]
\centering
\includegraphics{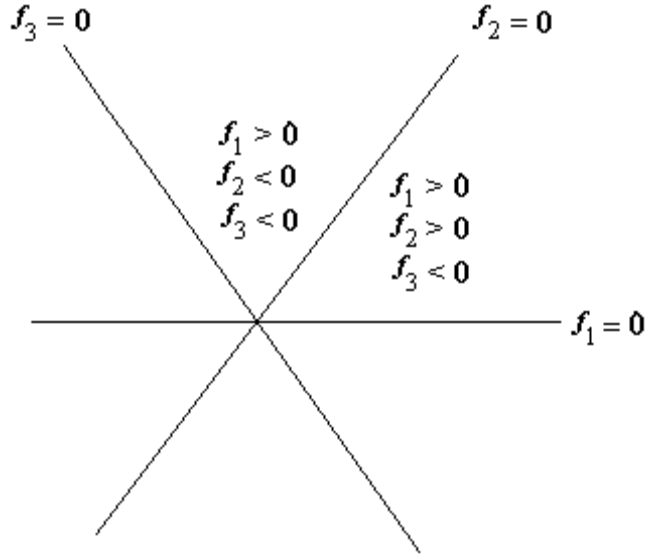}
\caption{\hspace{1mm} Here $\dim V = 2$ and $k = 3$. The three lines
divide the plane into $13$ cells ($6$ open sectors, $6$ open rays,
and the cell consisting of the origin. }
\end{figure}

1) A {\it cell} in $V$ with respect to $\mathbb{H}$ is a non-empty
set $A$ obtained by choosing, for each $i$, the half-space $f_i > 0$
or $f_i < 0$ or the hyperplane $f_i = 0$.

2) The cells defined by taking only the $f_i$ corresponding to
half-spaces are called {\it chambers}. Essentially, the chambers can
be described as the cells of maximal dimension. Sometimes what we
defined here as chambers are called cells in the literature.

3) Chambers are the connected components  of the complement
\begin{equation}
    V \setminus \bigcup\limits_{i=0}^k{H}_i.
\end{equation}

4) Let $C$ be the simplicial cone in $V$ defined by the inequalities
\begin{equation}
    f_i \geq 0 \text{ for } i = 1,2,\dots,n.
\end{equation}
It is called the {\it fundamental chamber}.

\ssec{The Coxeter complex}\label{Coxeter_complex} A cell $B$ is said
to be a {\it face} of $A$ if its description is obtained from that
of $A$ by replacing several inequalities by equalities. In this
case, we write
\begin{equation}
    B \leq A
\end{equation}
and this relation is said to be {\it face relation}. We have
\begin{equation*}
    \overline{A} = \bigcup\limits_{B \leq A}B,
\end{equation*}
and
\begin{equation}
  \label{cells_incl}
    B \leq A \Longleftrightarrow   \overline{B} \subseteq \overline{A}
\end{equation}

Let $\Sigma$ be the poset consisting of the open cells, ordered  by
the face relation. By (\ref{cells_incl}) $\Sigma$ is isomorphic to
the set of {\it closed} cells, see \cite[Ch.I]{Br}.

\sssec{The simplicial complex $\Sigma(W,S)$} Let $(W, S)$ be a
Coxeter group. Consider the subcomplex $\Sigma_{\leq C}$ consisting
of the faces of $C$. With every face $A \leq C$, we associate its
stabilizer $W_A = \{w \in W\mid wA = A\}$. By a theorem  in
\cite[Ch.I, \S5F]{Br}, $W_A$ is generated by a subset  $A\subset S$.

There is a function $\Phi$ from $\Sigma_{\leq C}$ to the set of
special subgroups of $W$; this $\Phi$ is a bijection (\cite[Ch.I,
\S5H]{Br}):
\begin{equation}
  \label{simplex_iso}
    \Sigma_{\leq C} \approx \text{(special subgroups of $(W,S)$)}^{op},
\end{equation}
where \lq\lq op" indicates that we are using the opposite of the
usual order on the set of special subgroups.

\begin{figure}[h]
\centering
\includegraphics{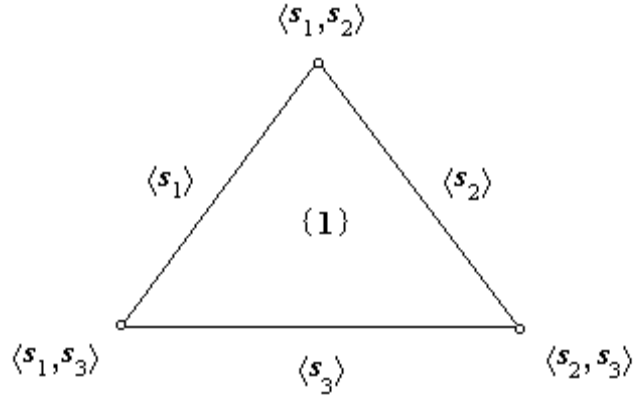}
\caption{\hspace{1mm} Illustration of the isomorphism
(\ref{simplex_iso}) for $n = 3$.}
\end{figure}

The $W$-action can be used to extend the isomorphism
(\ref{simplex_iso}) to an isomorphism of the whole poset $\Sigma$
with the set of {\it special cosets} in $W$, i.e., the cosets $wW'$
of special subgroups.

\begin{Theorem}{\em(\cite[Ch.I, \S5H]{Br})}
  \label{th_special_cosets_iso}
There is a poset isomorphism
\begin{equation}
  \label{special_cosets_iso}
    \Sigma \approx \text{(special cosets of $(W,S)$)}^{op}
\end{equation}
compatible with the $W$-action on the special cosets by
left-translation. 
\end{Theorem}

The Theorem \ref{th_special_cosets_iso} allows one to introduce
geometry into abstract group theory. Let $W$ be a group, possibly
infinite, generated by a subset $S$ consisting of elements of order
2. Define, a {\it special coset}  to be a coset
$w\langle{S'}\rangle$ with $w \in W$ and $S' \subset S$. Now define
$\Sigma = \Sigma(W, S)$ to be the poset of special cosets, ordered
by the opposite of the inclusion relation: $B \leq A$ in $\Sigma$ if
and only if $B \supset A$ as subsets of $W$, in which case we say
that $B$ is a {\it face} of $A$.

Following Tits, $\Sigma$ is called the {\it\bf Coxeter complex}
associated to $(W, S)$, it is also called the \lq\lq {\it apartment}
associated to $(W, S)$", see \cite{Br02}.

\ssec{The Euler characteristic}
   \label{euler}
For any finite simplicial complex, the {\it Euler characteristic} is
defined as the alternating sum
\begin{equation*}
    \chi = k_0 - k_1 + k_2 - k_3 + ...,
\end{equation*}
where $k_n$ denotes the number of cells of dimension $n$ in the
complex.

For any topological space, we can define the $n$th {\it Betti
number} $b_n$ as the rank of the $n$-th singular homology group. The
Euler characteristic is then equal to the alternating sum
\begin{equation*}
    \chi = b_0 - b_1 + b_2 - b_3 + ....
\end{equation*}
This quantity is well-defined if the Betti numbers are all finite
and if they are zero beyond a certain index $n_0$.

\section{The Poincar\'e series of the Lann\'er and quasi-Lann\'er
groups (new results)} Having computed something different from
Worthington's results,  we realized that means of verification are
badly needed. Besides, later we obtained a bit different picture
describing distribution of poles than the one Gal gave for $QL10_2$.
But our goal was not to refute (or verify) somebody's results but to
say something new. At first, we could only say something negative
(\lq\lq there is no reciprocity", \lq\lq not all non-real poles lie
on the unit circle $C$ centered at the origin", etc.), which was not
appealing. Fortunately, we managed to observe several patterns that
one can formulate in a positive way:
\begin{enumerate}
\item {\sl If the number of vertices of a given quasi-Lann\'er diagram is even,
the Euler characteristic vanishes}.
\item {\sl The difference of degrees of the numerator and denominator of the
 Poincar\'e series
 is always $\leq 1$ in the quasi-Lann\'er cases}.
\item {\sl The virgin form\footnote{The virgin form is defined in
subsection \ref{vf}; the reduced form of the rational Poincar\'e
series is its representation as an irreducible fraction.} of the
Poincar\'e series is equal to its reduced form in the quasi-Lann\'er
cases bar the following exceptions: $QL4_8$, $QL4_{12}$, $QL4_{14}$,
$QL4_{15}$, $QL4_{19}$, $QL6_5$, $QL6_9$, $QL6_{11}$, $QL8_1$,
$QL8_2$, $QL8_4$}.
\end{enumerate}
 In what follows we give {\it a priori} proofs of these and several
 other patterns.

\subsection{Reciprocity for the Lann\'er diagrams}

\sssbegin{Lemma}[On reciprocity]
 \label{lem_recipr}
 Let a polynomial $S(t)$ be factorized as follows:
\[
 S(t) = U(t)V(t),
\]
and let $U(t)$ be anti-reciprocal. Then
\begin{equation*}
 \begin{array}{cc}
 S(t) \text{ is anti-reciprocal if and only if } V(t) \text{
 is reciprocal}, \\
 S(t) \text{ is reciprocal if and only if } V(t) \text{
 is anti-reciprocal}.
 \end{array}
\end{equation*}
\end{Lemma}
Since $t - 1$ is an anti-reciprocal polynomial, this rather obvious
lemma helps us to understand that denominators of Poincar\'e series
of all Lann\'er diagrams on $4$ vertices are anti-reciprocal, see
Tables \ref{lanner_41_44} and \ref{lanner_46_49}.

\sssbegin{Proposition}\label{arecip} $1)$ Let $(W, S)$ be the
Coxeter system, such that all special subgroups of $W$ are finite.
If $\card{S}$ is even, the Poincar\'e series $W(t)$ is {\bf
anti-reciprocal}. If $\card{S}$ is odd, the Poincar\'e series $W(t)$
is {\bf reciprocal}.

$2)$ The Poincar\'e series $W(t)$ of the Lann\'er groups on $4$
vertices are anti-reciprocal, and the Poincar\'e series $W(t)$ of
the Lann\'er groups on $5$ vertices are reciprocal.

\end{Proposition}

\begin{proof}
1) By \eqref{exBu} and \eqref{exSte} we have
\begin{equation}
 \label{antirecpr_2}
 \begin{array}{ccc}
 & -\displaystyle\frac{\epsilon(D)}{W_D(t)} =& \sum\limits_{X\subsetneq D}
 \displaystyle\frac{\eps(X)}{W_X(t)}, \\
 \\
 & \displaystyle\frac{1}{W_D(t^{-1})} =&
 \sum\limits_{|W_X|<\infty}\displaystyle\frac{\eps(X)}{W_X(t)}.\\
 \end{array}
\end{equation}
 Since all special subgroups are finite, then the right hand sides in both
 equations coincide. Thus, if $|S|$ is even (resp. odd), then $-\epsilon(D)$ is negative
 (resp. positive) and the Poincar\'e series is anti-reciprocal (resp. reciprocal).

 2) Note, that in the case of Lann\'er groups,
 \begin{equation}
 \label{note_lanner}
\begin{minipage}[l]{13cm}
the set of all finite special subgroups coincides with the set of
all special subgroups.
\end{minipage}\end{equation}
 Then the statement desired follows from heading 1).
 \end{proof}

\sssbegin{Conjecture}
 \label{conj2}
For the Coxeter groups,  the (anti)\-reciprocity never holds, bar
the cases listed in Prop. $\ref{arecip}$.\end{Conjecture}

\subsection{When does the Euler characteristics vanish?}

\sssbegin{Proposition}
   \label{prop_Euler}
 The Euler characteristics $\chi(W)$ of the group $W$ vanishes (equivalently,
the denominator of the Poincar\'e series has the root $t = 1$)
 in the following cases:

 $1)$ For any affine Coxeter group.

 $2)$ For any infinite
 (non-affine) Coxeter group $(W, S)$ with $|S|$ even.
(Of course, this case includes (quasi-)Lann\'er groups whose Coxeter
diagrams have
 even number of vertices.)
\end{Proposition}

 \begin{proof}
 1) Follows from \eqref{HPWloopquot}.

 2) By \eqref{antirecpr_2}, we have
\begin{equation*}
 -\frac{\epsilon(D)}{W_D(t)} = \frac{1}{W_D(t^{-1})}.
\end{equation*}
 Substituting $t = 1$ for $\card{D}$ even, we see that
\begin{equation*}
 -\frac{1}{W_D(1)} = \frac{1}{W_D(1)}, \quad \text{ i.e.,} \quad \frac{1}{W_D(1)} = 0.
\end{equation*}
For illustration of this fact, see Tables \ref{qlanner_ql41_ql49},
 \ref{qlanner_ql410_ql416}, \ref{qlanner_ql417_ql423},
 \ref{q_lanner_61}
 and \ref{q_lanner_67_69}.
 \end{proof}


\sssec{The virgin form of the numerator}\label{vf} The numerator of
$W_D(t)$ is equal to the denominator of the sum
$\mathop{\sum}\limits_{X\subsetneq D} \frac{\eps(X)}{W_X(t)}$. By
\eqref{solo}, for the finite Coxeter group $W_X$ with exponents
\begin{equation*}
 m_1, m_2, \dots, m_k,
\end{equation*}
 the Poincar\'e series $W_X$ is a polynomial of the form
\begin{equation}
 \label{[n]s}
[m_1 + 1][m_2 + 1]\dots [m_k + 1].
\end{equation}
The least common multiple
\begin{equation}\label{virg}{\rm Virg}(D):=\mathop{{\rm LCM}}\limits_{X\subsetneq
D\text{ such that } |W_X|<\infty}W_X(t)\end{equation} is said to be
the {\it virgin form} of the numerator of $W_D(t)$.

\parbegin{Lemma} The Poincar\'e series $W_D(t)$ can be expressed as a rational fraction
whose numerator is ${\rm Virg}(D)$. \end{Lemma}

\begin{proof} The statement is obvious if all special subgroups $W_X$ are finite:
then the numerator of $W_D(t)$ is equal to the denominator of the
sum $\mathop{\sum}\limits_{X\subsetneq D} \frac{\eps(X)}{W_X(t)}$
and all denominators of its summands are polynomials of the form
\eqref{[n]}. The general case is done by induction on $|X|$.
\end{proof}

\parbegin{Corollary} Let $\displaystyle\frac{\eps(X)}{W_X(t)}$ be expressed as an irreducible
fraction. Then the LCM of all denominators in the sum
$\mathop{\sum}\limits_{X\subsetneq D}
\displaystyle\frac{\eps(X)}{W_X(t)}$ is equal to ${\rm Virg}(D)$.
\end{Corollary}

\begin{proof} Indeed, if $|W_X|=\infty$, then the denominator of
the irreducible fraction $\displaystyle\frac{\eps(X)}{W_X(t)}$
divides ${\rm Virg}(X)$ and ${\rm Virg}(X)$ divides ${\rm Virg}(D)$.
If $|W_X|<\infty$, then $W_X(t)$ divides ${\rm Virg}(D)$ by
definition. Hence, the LCM of denominators divides ${\rm Virg}(D)$.

Implication in the opposite direction: divisibility of the LCM of
denominators by ${\rm Virg}(D)$ is obvious.
\end{proof}

If $|W_X|<\infty$, then $W_X(t)$ is of the form \eqref{[n]}. We
would like to represent ${\rm Virg}(D)$ in the same form, but this
is not always possible: if $m$ and $n$ are not relatively prime,
then $[m]$ and $[n]$ are not relatively prime. On the other hand,
each polynomial $[n]$ can be represented as the product of
irreducible over $\Qee$ cyclotomic polynomials $\Phi_i(t)$, where
$i=2, 3,\dots$, namely
\begin{equation}\label{virg1}[n]=\mathop{\prod}\limits_{i|n, ~i>1}
\Phi_i(t).\end{equation}
\begin{center}
\begin{minipage}[c]{12cm}
\begin{center}THEREFORE, IT IS NATURAL TO COMPUTE ${\rm Virg}(D)$ IN THE FORM OF
THE PRODUCT OF THE $\Phi_i(t)$. \end{center}
\end{minipage}
\end{center}
It is convenient to introduce one
more notation:
\begin{equation}
 \label{[n]'}
[n']:=1+t^n; \text{~~ observe that~}[n][n']=[2n].
\end{equation}

\ssec{Degrees of the denominators} In this section, we define the
polynomials $P,Q,R,S$ by setting:
\begin{equation}
 \label{WiWinv}
W(t):=\frac{R(t)}{S(t)}\text{~~and ~~}W(t^{-1}):=\frac{P(t)}{Q(t)}.
\end{equation}

\sssbegin{Proposition}
 \label{prop_Meb1} For any Coxeter group, we have

 $1)$ $\deg P=\deg Q$;

 $2)$ $\deg S < \deg R$ if and only if $t\mid Q(t)$.
\end{Proposition}

\begin{proof}
 1) According to the Solomon formula \eqref{solo},
 for any finite Coxeter group, every cyclotomic polynomial-factor of $W(t^{-1})$
 turns into the fraction
\[
  \frac{1 + t + ... + t^{n-1}}{t^{n-1}}
\]
For any affine Coxeter group, Proposition follows from
\eqref{HPWloopquot}. For any infinite Coxeter group, we use the
Steinberg formula \eqref{exSte}. Recall that $1$ in the numerator
above is the summand corresponding to the empty set in
\eqref{exSte}. This summand contributes the maximal degree equal to
the degree of the denominator of $W(t^{-1})$. Therefore, $\deg
P=\deg Q$.

2) Note that $t\mid Q(t)$ if and only if $a_0 = 0$, where $a_0$ is
the constant term of $Q(t)$, which becomes the highest coefficient
of $S(t)$ under the substitution $t \longmapsto t^{-1}$. Thus, the
condition $t\mid Q(t)$ is equivalent to $\deg S < \deg R$.
\end{proof}

\sssbegin{Proposition}
 $1)$ For the function $\epsilon(X) = (-1)^{|X|}$, we have
 \begin{equation}
 \label{eps_func1}
 a)  \mathop{\sum}\limits_{\emptyset \subseteq X \subseteq D}\epsilon(X) = 0,
  \qquad
 b) \mathop{\sum}\limits_{\emptyset \subseteq X \subseteq D}\epsilon(X)\,|X| = 0.
 \end{equation}

$2)$ We have:
 \begin{equation}
 \label{prop_Moeb_func2}
 \displaystyle\frac{1}{W(t^{-1})|_{t=0}}=
 \mathop{\sum}\limits_{X \subset D ~\mid~
 |W_X| < \infty}\epsilon(X).
\end{equation}

$3)$ If $f(t)$ is the product of $k$ factors $[n_i]$, then
$f'(0)=k$.
 \begin{equation}
  \label{deriv_1}
 f(t) = \prod\limits_{i=1}^{k}[n_i] \quad \Longrightarrow \quad f'(0) = k.
 \end{equation}

 $4)$ If $|X|<\infty$,  then
 \begin{equation}
  \label{deriv_2}
 W'_X(t) \mid_{t=0} =  |X|.
 \end{equation}

\end{Proposition}
\begin{proof}
1) Formula (\ref{eps_func1},a) holds since
\begin{equation*}
 \sum\limits_{\emptyset \subseteq X \subseteq D}\epsilon(X) =
    \sum\limits_{k = 0}^{n}(-1)^k \binom{n}{k} = (1 - 1)^n = 0,
 \end{equation*}
 and formula (\ref{eps_func1},b) is true since
 \begin{equation*}
 \sum\limits_{\emptyset \subseteq X \subseteq D}\epsilon(X)\mid X \mid =
    \sum\limits_{k = 0}^{n}(-1)^k \binom{n}{k}k = -n\sum\limits_{k = 0}^{n}
    (-1)^{k-1} \binom{n-1}{k-1}
    = -n(1 - 1)^n = 0.
 \end{equation*}

 2) By the Steinberg formula \eqref{exSte} for
$\displaystyle\frac{1}{W(t^{-1})}$, every summand of \eqref{exSte}
is of the form
 \begin{equation}
 u_i(Z) = \frac{\eps(Z)}{[n_1][n_2]\dots[n_k]}
 \end{equation}
 Since $[n_i]|_{t=0} = 1$, we have $u_i(Z) = \eps(Z)$.

 3) Eq. \eqref{deriv_1} holds since
\begin{equation*}
  f'(t) = \sum\limits_{i=1}^{k}[n_i]'\prod\limits_{j=1; j\neq i}^{k}[n_i], \quad \text{ and } \quad
  f'(0) = \sum\limits_{i=1}^{k}[n_i]'\mid_{t=0} =
\sum\limits_{i=1}^{k} 1  = k.
 \end{equation*}

 4) Since the number of factors $[n_i]$ in the Poincar\'e series
of any finite Coxeter group is equal to the number of its
generators, i.e., to the number of vertices $|X|$, then
$\eqref{deriv_2}$ follows from $\eqref{deriv_1}$.
\end{proof}

\sssbegin{Proposition}
 \label{th_Meb}
 $1)$ The following relation holds:
 \begin{equation}
 \label{eq_func_1}
 \displaystyle\frac{1}{W(t^{-1})}\Big |_{t=0}
 ~~=~~
 (-1)^{|D|+1}-\mathop{\sum}\limits_{|W_X|=\infty}(-1)^{|X|}.
 \end{equation}
$2)$ For degrees of the numerator and denominator of $W(t)$, we
have:
 \begin{equation}
 \label{eq_func_2}
 \deg S<\deg R \quad \text{ if and only if } \quad
 \mathop{\sum} \limits_{|W_X|=\infty} \eps(X) =
 \mathop{\sum} \limits_{|W_X|=\infty}(-1)^{|X|}=(-1)^{|D|+1}.
 \end{equation}
$3)$ We have:
 \begin{equation}
 \label{eq_func_3}
 \left ( \frac{1}{W(t^{-1})} \right )' \Big |_{t=0} =
   \sum\limits_{|W_X|=\infty}{(-1)}^{|X|}|X| + {(-1)}^{|D|}|D|.
 \end{equation}
\end{Proposition}

\begin{proof}
1) Follows from the fact that the sum in \eqref{prop_Moeb_func2}
differs from (\ref{eps_func1},a) by summands associated with the
infinite special subgroups and the subset $X = D$.

2) According to Proposition \ref{prop_Meb1} the condition $\deg
S<\deg R$ is equivalent to $Q(t) = 0$, or, in other words, to
$$
\displaystyle\frac{1}{W(t^{-1})}\Big |_{t=0} = 0.
$$
Then the statement follows from \eqref{eq_func_1}.

3)  Since
$$
   \frac{1}{W(t^{-1})} = \sum\limits_{|W_X| <
   \infty}\frac{\eps(X)}{W_X(t)},
$$
we have
\begin{equation*}
  \left ( \frac{1}{W(t^{-1})} \right )' = -\sum\limits_{|W_X| <
   \infty}\frac{W'_X(t)\eps(X)}{\left(W_X(t)\right)^2}
\end{equation*}
For any finite Coxeter group, we have $W_X(t)|_{t=0} = 1$. By
\eqref{deriv_2} $W'_X(t)\Big |_{t=0} = |X|$, and we see that
\begin{equation*}
\left ( \frac{1}{W(t^{-1})} \right )'\Big |_{t=0}=
-\sum\limits_{|W_X| < \infty}\eps(X)|X|.
\end{equation*}
According to (\ref{eps_func1},b) we have:
\begin{equation*}
\left ( \frac{1}{W(t^{-1})} \right )'\Big |_{t=0} =
\sum\limits_{|W_X| =
   \infty}\eps(X)|X| + \eps(D)|D|,
\end{equation*}
and \eqref{eq_func_3} holds.

 \end{proof}

\sssbegin{Corollary}
 \label{Meb5} For the quasi-Lann\'er
diagrams, we have $\deg S<\deg R$ only in the following cases:
 \begin{equation}\label{S<R}
\begin{array}{l}
QL4_1, \quad QL4_2, \quad QL4_3, \quad QL4_6, \quad QL4_8, \quad
 QL4_{14}, \quad QL4_{17}, \quad QL4_{18}, \quad QL4_{19};\\
QL5_1, \quad QL5_2, \quad QL5_3, \quad QL5_6, \quad QL5_7;\\
QL6_5, \quad QL6_7, \quad QL6_{10};\\
QL7_1, \quad QL7_2, \quad QL7_3; \\
QL8_1, \quad QL8_2, \quad QL8_3, \quad QL8_4;\\
QL9_1, \quad QL9_2, \quad QL9_3;\\
QL10_1.\\
\end{array}
 \end{equation}
\end{Corollary}

 \begin{proof}
Indeed, only in these cases there is a {\bf single infinite special
subgroup} in the given quasi-Lann\'er Coxeter group.

For diagrams on $4$ vertices, this subgroup is associated with $X$,
such that $|X| = 3$. Further, $|D| + 1 = 5$, and $(-1)^{|X|} =
(-1)^{|D| + 1} = -1$. Then, the statement follows from
\eqref{eq_func_2}. The cases of $>4$ vertices are absolutely
analogous.
 \end{proof}

\ssec{The coefficients $b_n$ and $b_{n-1}$ of the denominator}
  Let $b_n$ (resp. $b_{n-1}$) be the coefficient corresponding
  the degree $n$ (resp. $n-1$) of the denominator of Poincar\'e
  series $W(t)$. Consider
\begin{equation}
 \label{def_ai}
     \frac{1}{W(t^{-1})} = \frac{Q(t)}{P(t)}, \text{ where }
     Q(t) = \sum\limits_{i=0}^n{a}_i{t}^i, \quad P(t) =
    \prod\limits_{i = 1}^{m}[n_i].
 \end{equation}
We have:
\begin{equation}
    \label{bi_ai}
    a_0 = b_n, \qquad a_1 = b_{n-1}.
 \end{equation}
Now, we will prove two theorems predicting values of coefficients
$b_n$ and $b_{n-1}$ of the denominator. It is clear  that the other
coefficients $b_i$ of the denominator can be predicted in the same
way. Note that calculations of $b_n$ and $b_{n-1}$ are closely
connected with the poset of infinite special subgroups in $W$. The
following theorem is devoted to the coefficient $b_n$. Actually, the
conclusion \eqref{eq_func_2} is a particular case of this theorem.

\sssbegin{Theorem}
   \label{th_Highest}
All Poincar\'e series in the tables are normalized so that $b_n =
1$.

$1)$ For the coefficient $b_n$ of the highest term of the
denominator $S(t)$ of $W(t)$, we have
\begin{equation}
 \label{eq_highest}
   b_n = (-1)^{|D| + 1} - \sum\limits_{|W_X| = \infty}(-1)^{|X|}\quad.
 \end{equation}

$2)$ For any $0$-terminal Coxeter group, in particular, for any
Lann\'er group, we have
\begin{equation}
     \label{eq_highest_Lanner}
   b_n = (-1)^{|D| + 1},
 \end{equation}
see Tables $\ref{lanner_41_44}$, $\ref{lanner_46_49}$,
$\ref{lanner_51}$.

$3)$ For any $1$-terminal Coxeter group $W$, in particular, for any
quasi-Lann\'er group, we have
\begin{equation}
     \label{eq_highest_qLanner}
   b_n = (-1)^{|D|}(\Infi - 1), \quad \Infi = b_n + 1,
 \end{equation}
where $\Infi$ is the number of infinite special subgroups in $W$,
see Tables $\ref{qlanner_ql41_ql49}$ -- $\ref{q_lanner_9_10a}$.

$4)$ For any $k$-terminal  Coxeter group $W$, we have
\begin{equation}
     \label{eq_highest_qLanner_common}
   b_n = (-1)^{|D|+1}\sum\limits_{i=0}^m{(-1)^{m}\Infi_m},
 \end{equation}
\end{Theorem}

\begin{proof}
1) Recall that by \eqref{WiWinv} $P(t)$ (resp. $Q(t)$) is the
numerator (resp. denominator) of $W(t^{-1})$. The case $P(t) = 0$ is
considered in Theorem \ref{th_Meb}. Now, let $P(t) \neq 0$. Since
\begin{equation*}
Q(t) = \prod\limits_{i = 1}^k[n_i],
\end{equation*}
we have $Q(0) = 1$. Thus,
\begin{equation*}
 \frac{1}{W(t^{-1})}\Big |_{t=0} = P(0) = a_0 \neq 0,
 \end{equation*}
where $a_0$ is the constant term of the denominator $P(t)$ of
$W(t^{-1})$. Substitution $t \longmapsto t^{-1}$ turns  $a_0$ into
$b_n$, the coefficient of the highest term of the denominator $R(t)$
of $W(t)$.

2) For Lann\'er groups (and also $0$-terminal) the term
$\sum\limits_{|W_X| = \infty}(-1)^{|X|}$ in \eqref{eq_highest}
vanishes.

3) For quasi-Lann\'er (and $1$-terminal) groups  we have
 $|D| = |X| + 1$, where $X$ is the subdiagram corresponding any infinite
subgroups, and \eqref{eq_highest_qLanner} holds.
 \end{proof}

\subsubsection{The $[n]$-complete and reduced forms of the Poincar\'e
series}
  \label{n-comp}
  The following theorem is devoted to predicting
the coefficient $b_{n-1}$ of the denominator of the Poincar\'e
series. The calculation of $b_{n-1}$ is based on the parameter $m$
of the numerator meaning the number of factors like $[n_i]$ in the
numerator. However, the numerator which we consider is not mandatory
irreducible. If it contains a divisor of some $[n_i]$, we multiply
the numerator and denominator to get only factors like $[n_i]$. This
non-irreducible form of the the Poincar\'e series is said to be the
{\bf $[n]$-complete} form. Thus, our prediction is related to the
numerator of the $[n]$-complete form. Note that

(1) For the quasi-Lann\'er Coxeter groups, there are only two cases,
namely $QL8_1$ and $QL8_2$,  with two $[n]$-incomplete factors. In
the cases $QL4_8$, $QL4_{12}$, $QL4_{14}$, $QL4_{15}$, $QL4_{19}$,
$QL6_5$, $QL6_9$, $QL6_{11}$, $QL8_4$, there is only one
$[n]$-incomplete  factor. In the remaining cases the $[n]$-complete
form and reduced form coincide.

(2) The following fact holds: after reduction of the $[n]$-complete
form of the Poincar\'e series for Lann\'er and quasi-Lann\'er groups
the number of factors $m$ is not changed. None of the factors
$[n_i]$ is completely reduced.

(3) The difference of degrees of the numerator and denominator $\deg
R - \deg S$ does not change under  reduction. This fact allows us to
calculate $\deg R - \deg S$ for the $[n]$-complete form and to apply
it to the reduced form.

\sssbegin{Theorem}
   \label{th_b_n-1}
Let $m$ be the number of factors $[n_i]$ in the $[n]$-complete form
as \eqref{def_ai}. Then
\begin{equation}
 \label{eq_2_coef}
 b_{n-1} - m b_n  = \sum\limits_{|W_X| = \infty}(-1)^{|X|}|X| + (-1)^{|D|}|D|.
 \end{equation}

\end{Theorem}

\begin{proof}
   Since
\begin{equation*}
     \left ( \frac{1}{W(t^{-1})} \right )'\Big |_{t=0}  \quad=\quad
     \frac{f'(t)g(t) - f(t)g'(t)}{(g(t))^2} \Big |_{t=0},
     \text{ and } g(0) = 1,
 \end{equation*}
 and by \eqref{deriv_1} $g'(0) = m$, we have
\begin{equation*}
     \left ( \frac{1}{W(t^{-1})} \right )'\Big |_{t=0}  \quad=\quad
      f'(0) - f(0)g'(0) \quad=\quad a_1 - m{a}_0.
 \end{equation*}
By \eqref{bi_ai},
\begin{equation*}
     \left ( \frac{1}{W(t^{-1})} \right )'\Big |_{t=0}   \quad=\quad b_{n-1} - m{b}_n.
 \end{equation*}
 Then eq. \eqref{eq_2_coef} follows from \eqref{eq_func_3}.
\end{proof}

\sssbegin{Corollary}
   For any Lann\'er group we have:
\begin{equation}
  \label{col_2_coef}
   \begin{split}
      & b_n = -1, \quad b_{n-1} = 1 \quad \text{ for } \quad L4_i, \quad 1 \leq i \leq 9, \\
      & b_n = 1, \quad b_{n-1} = -1  \quad \text{ for } \quad L5_i, \quad  i = 1,3,4, \\
      & b_n = 1, \quad b_{n-1} = 0  \quad \text{ for } \quad L5_i, \quad  i = 2,5.
   \end{split}
 \end{equation}
\end{Corollary}

\begin{Remark} Eq. \eqref{col_2_coef} holds for the $[n]$-complete form, and
does not hold for the reduced form, see Tables \ref{lanner_41_44},
\ref{lanner_46_49}, \ref{lanner_51}. All Poincar\'e series in Tables
\ref{lanner_41_44}, \ref{lanner_46_49}, \ref{lanner_51} are
normalized so that $b_n = 1$. \end{Remark}

\begin{proof}
   Since Lann\'er groups does not contain infinite subgroups, i.e., $\Infi = 0$,
   then $b_n = (-1)^{|D| + 1}$, see \eqref{eq_highest_Lanner}.

   For $|D| = 4$, the number of factors $m = 3$. In this case, by
   \eqref{eq_2_coef}  we have
\begin{equation*}
      b_{n-1} = m b_n + (-1)^{|D|}|D| = 3(-1) +  4 = 1.
 \end{equation*}

   For $|D| = 5$, the number of factors $m = 4$ (except for $L5_2$ and $L5_5$).
   In this case, by \eqref{eq_2_coef}  we have
\begin{equation*}
      b_{n-1} = m b_n + (-1)^{|D|}|D| = 4 -  5 = 1.
 \end{equation*}

For $|D| = 5$, cases $L5_2$ and $L5_5$, the number of factors $m =
5$. In this case, by \eqref{eq_2_coef}  we have
\begin{equation*}
      b_{n-1} = m b_n + (-1)^{|D|}|D| = 5 -  5 = 0.
 \end{equation*}

 \end{proof}

\sssbegin{Corollary}
 \label{cor_sindle_inf_gr}
$1)$ For any Coxeter group with a single infinite subgroup, we have:
\begin{equation}
 \label{bn_and_b_n-1}
 b_n = 0, \quad \text{ and } \quad b_{n-1} \neq 0.
\end{equation}
In this case,  $\deg R - \deg S=1$.

$2)$ For any quasi-Lann\'er group (and also for any $1$-terminal
Coxeter group), the difference of degrees of the numerator and
denominator of the Poincar\'e series is $\deg R - \deg S\leq 1$.
\end{Corollary}

\begin{proof}
1) According to \eqref{eq_2_coef}, and since $\Infi = 1$, we have
\begin{equation*}
 b_{n-1} - m b_n  = (-1)^{|X|}|X| + (-1)^{|D|}|D|,
\end{equation*}
 where $|X| = |D| - 1$, i.e.,
\begin{equation*}
 b_{n-1} - m b_n  = (-1)^{|X|}(|X| - |X| + 1|) = (-1)^{|X|},
\end{equation*}
From \eqref{eq_highest_qLanner} we have $b_n = 0$, and therefore
$b_{n-1} = (-1)^{|X|}$, so \eqref{bn_and_b_n-1} holds.

2) Let $b_n = 0$. Since in the case of $1$-terminal Coxeter group
$|X| + 1 = |D|$ for all infinite subgroups $X \subset D$, all
summands ${(-1)}^{|X|}$  in \eqref{eq_highest} have the same sign.
Since $b_n = 0$, there exists only one infinite special subgroup,
and by 1) we have $\deg R - \deg S = 1$.
\end{proof}

\sssbegin{Conjecture} For ANY infinite Coxeter group, $\deg R - \deg
S\leq 1$.
\end{Conjecture}

\sssbegin{Corollary} Any infinite Coxeter group having exactly two
infinite special subgroups is $1$-terminal or $2$-terminal.

$1)$ For any $1$-terminal Coxeter group with exactly two infinite
special subgroups, we have:
\begin{equation}
 \label{2_Inf_gr}
 b_n = {(-1)}^{|D|}, \quad \text{ and } \quad b_{n-1} = {(-1)}^{|D|}(m + 2 - |D|).
\end{equation}

$2)$  The quasi-Lann\'er groups with exactly two infinite special
subgroups are as follows:

a)  For $|D| = 4$, $m = 3$, we have $b_n = 1, b_{n-1} = 1$, (see
cases $QL4_4$, $QL4_7$, $QL4_{11}$, $QL4_{12}$).

b) For $|D| = 4$, $m = 2$, we have $b_n = 1, b_{n-1} = 0$, (see
cases $QL4_5$, $QL4_{10}$, $QL4_{21}$).

c) For $|D| = 5$, $m = 4$, we have $b_n = -1, b_{n-1} = -1$, (see
case $QL5_8)$.

d) For $|D| = 6$, $m = 5$, we have $b_n = 1, b_{n-1} = 1$, (see
cases $QL6_1$, $QL6_3$, $QL6_{11}$).

e) For $|D| = 9$, $m = 8$, we have $b_n = -1, b_{n-1} = -1$, (see
case $QL9_4$).

f) For $|D| = 10$, $m = 9$, we have $b_n = 1, b_{n-1} = 1$, (see
case $QL10_2$).

$3)$ For any $2$-terminal Coxeter group with exactly two infinite
special subgroups, we have
\begin{equation}
 \label{2_Inf_gr_2-term}
 b_n = {(-1)}^{|D|+1}, \quad \text{ and } \quad b_{n-1} = {(-1)}^{|D|+1}(m + 1 - |D|).
\end{equation}
\end{Corollary}

\begin{proof} Follows from \eqref{eq_highest} and \eqref{eq_2_coef}.
\end{proof}


\ssec{Nerves and geometric realization of the group}\label{Nerve}
From \cite[p.374]{ChD}: The {\it nerve} of $(W, S)$, denoted by $N$,
is the poset of subsets $X\subset S$ for which the group $W_{X}:=(W,
X)$ is finite. The poset is ordered with respect inclusion. The {\it
proper nerve} of $(W, S)$ is the poset $N_{> \emptyset}$ consisting
of the nonempty subsets $X\subset S$ such that $W_{X}$ is finite.

Clearly, $N_{> \emptyset}$ is a {\it simplicial complex}. More
precisely, it is isomorphic to the poset of simplices of a
simplicial complex with vertex set $S$. (For more facts and
explanations, see subsect. \ref{Coxeter_complex}.)

For any finite poset $K$, let $\chi(K)$ denote the Euler
characteristic of its {\it geometric realization} (see subsect.
\ref{geom_realiz}, \ref{euler}). The following formula due to Serre
\cite{Se} connects the Poincar\'e series of a given Coxeter group
and its Euler characteristic:
\begin{equation}
 \label{Serre_frm}
\frac{1}{W(1)} = \chi(W) \qquad\text{({\bf The Serre Formula})}.
\end{equation}

\section{The code to compute the Poincar\'e series and means of control}

To compute the Poincar\'e series, we used the {\it
Mathematica}-based code \texttt{subg} due to D.~Chapovalov
\cite{DCh} and double-checked with a code due to R.~Stekolshchik.

\ssec{Code \texttt{subg}} We rewrite the expression \eqref{exBu} in
the following form
\begin{equation}
 \label{exBu2}
W_D(t)=\frac{-\eps(D)}{\mathop{\sum}\limits_{X\subsetneq D}
\frac{\eps(X)}{W_X(t)}}.
\end{equation}
This formula enables one to express the Poincar\'e series $W_D(t)$
in terms of the finite groups listed in Table 1. The corresponding
recursion was automatically generated by the code \texttt{subg}. The
format and notation (improving Coxeter symbols) are designed so that
each step can be easily verified by a human, and, on the other hand,
these intermediate results can be copied to {\it Mathematica} in
order to derive the final answer.

\ssec{Poles} Having found the Poincar\'e series we determined their
poles by means of {\it Ma\-the\-ma\-ti\-ca} and Molotkov verified
our findings with the help of the code \texttt{pari}, see \cite{rf}.

\subsection*{Acknowledgements.} We are thankful to
R.~Grigorchuk, P.~de la Harpe,  B.~Okun, O.~ Shwartsman, and \'E.~
Vinberg for useful comments. We are thankful to A.~Chapovalov and
D.~Chapovalov, and to V.~Molotkov for invaluable help in computing.

\newpage

\clearpage
\section{Tables}

\begin{table*}[h]
 \begin{center}
\caption{\hspace{1mm}The exponents, Coxeter number, and
 the maximal length of the elements in the spherical Coxeter groups
 with connected Coxeter diagram}
 \vspace{1mm}
 \renewcommand{\arraystretch}{1.7}
\begin{tabular}{||c|c|c|c|c||}
 \hline\hline
 Coxeter & Lie & exponents &
 maximal length & Coxeter \cr
 group & algebra & $m_i$ & $l(w_0)=\sum m_i$ & number $h$ \\
 \hline\hline
 $A_n$ & $\fsl(n+1)$ & $1, 2, 3, \dots, n$ &
 $\displaystyle\frac{n(n+1)}{2}$ & $n+1$\\
 \hline
 $B_n$ & $\fo (2n)$ for $n\geq 2$ & $1, 3, \dots,
 2n-1$ & $n^2$ & $2n$\\
 \hline
 $C_n$ & $\fsp(2n)$ for $n\geq 2$ & $1, 3, \dots,
 2n-1$ & $n^2$ & $2n$ \\
 \hline
 $D_n$ & $\fo(2n+1)$ & $1, 3, \dots, 2n-3; \;
 n-1$ & $n(n-1)$ & $2(n-1)$ \\
 \hline
 $G_2$ & $\fg_2$ & $1, 5$ & $6$ & $6$ \\
 \hline
 $F_4$ & $\ff_4$ & $1, 5, 7, 11$ & $24$ & $12$ \\
 \hline
 $E_6$ & $\fe_6$ & $1, 4, 5, 7, 8, 11$ & $36$ & $12$ \\
 \hline
 $E_7$ & $\fe_7$ & $1, 5, 7, 9, 11, 13, 17$ & $63$ & $18$ \\
 \hline
 $E_8$ & $\fe_8$ & $1, 7, 11, 13, 17, 19, 23, 29$ & $120$ & $30$ \\
 \hline
 $I_2^{(m)}$ for $m > 6$ & --- & $1, m-1$ & $m$ & $m$ \\
 or $m = 5$ & & & & \\
 \hline
 $H_3$ & --- & $1, 5, 9$ & $15$ & $10$ \\
 \hline
 $H_4$ & --- & $1, 11, 19, 29$ & $60$ & $30$ \\
 \hline\hline
\end{tabular}
 \label{table_expon}
 \end{center}
 {\bf Note}. The groups $I_2^{(m)}$ are the non-crystallographic dihedral groups for
 $m = 5$ and $m > 6$. For $m = 3, 4, \text{ and } 6$, respectively,
 we have the crystallographic dihedral group as follows:
 \begin{equation*}
 A_2=I_2^{(3)}, \qquad
 B_2 = C_2=I_2^{(4)}, \qquad
 G_2=I_2^{(6)}.
 \end{equation*}
\end{table*}

 \begin{table*}[h]
 \begin{center}
\caption{\hspace{1mm}The Poincar\'e series of the spherical Coxeter
groups
 with connected Coxeter diagram}
 \vspace{1mm}
 \renewcommand{\arraystretch}{1.7}
\begin{tabular}{||c|c||}
 \hline\hline
 Coxeter group& its Poincar\'e series\cr
 \hline\hline
 $A_n$ & $[2]\dots[n+1]$\\
 \hline
 $B_n$ & $[2]\dots[2n]$\\
 \hline
 $D_n$ & $[2]\dots[2n-2][n]$ \\
 \hline
 $G_2$ & $[2][6]$ \\
 \hline
 $F_4$ & $[2][6][8][12]$ \\
 \hline
 $E_6$ & $[2][5][6][8][9][12]$\\
 \hline
 $E_7$ &$[2][6][8][10][12][14][18]$\\
 \hline
 $E_8$ & $[2][8][12][14][18][20][24][30]$\\
 \hline
 $I_2^{(m)}$ for $5\leq m \leq\infty$ &$[2][m]$\\
 \hline
 $H_3$ &$[2][6][10]$\\
 \hline
 $H_4$ & $[2][12][20][30]$\\
\hline\hline
 Coxeter group& multiparameter Poincar\'e series ($u:=t_1t_2$)\cr
 \hline\hline
$I_2^{(m)}$ for $6\leq m=2k \leq \infty$ &$(1+t_1)(1+t_2)
\begin{cases}\frac{1-u^m}{1-u}&\text{if $m<\infty$}\\
\frac{1}{1-u}&\text{if $m=\infty$}\end{cases}$\\
 \hline
 $F_4$ &
 $\left(\mathop{\prod}\limits_{j=1,2}(1+t_j)(1+t_j+t_j^2)(1+ut_j)\right)
 \mathop{\prod}\limits_{i=1}^3(1+u^i)$ \\
 \hline
 $B_n$ ($t_2$ corr. to the short root)& $\mathop{\prod}\limits_{i=0}^{n-1}(1+t_1^it_2)(1+t_1+\dots +t_1^i)$\\
\hline
$C_n$ ($t_2$ corr. to the long root)& $\mathop{\prod}\limits_{i=0}^{n-1}(1+t_1^it_2)(1+t_1+\dots +t_1^i)$\\
\hline\hline
\end{tabular}
 \label{table_grfunell}
 \end{center}
 \end{table*}


\begin{table*}[h]
 \begin{center}
\caption{\hspace{1mm}The Lann\'er diagrams on $4$ vertices and
Poincar\'e series}
 \vspace{1mm}
\begin{tabular}{||c|c|c||}
 \hline\hline
 & & \cr
 Label & Diagram & Poincar\'e series \cr
       & Degrees & $\chi=0$ in all cases \\
 & & \\
 \hline\hline
 & & \cr
 $L4_1$ & \includegraphics{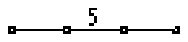} &
 $\displaystyle\frac{[2][6][10]}{t^{15} - t^{14} - t^{11} + t^4 + t - 1} = $\cr
 & $(11, 11)$ & \cr
 & & \footnotesize $\displaystyle\frac{[2][6][5']}{t^{11} - 2t^{10} + t^9 - t^7 + 2t^6 - 2t^5 + t^4 - t^2 + 2t - 1}$ \cr
\hline\hline
\end{tabular}
 \label{lanner_41_44}
 \end{center}
 \vspace{1mm}
 \end{table*}

\newpage

\begin{table*}[h]
 \begin{center}
\caption{\hspace{1mm}The Lann\'er diagrams on $4$ vertices and
Poincar\'e series (cont.)}
 \vspace{1mm}
\begin{tabular}{||c|c|c||}
 \hline\hline
 & & \cr
 Label & Diagram & Poincar\'e series \cr
       & Degrees & $\chi=0$ in all cases \\
 & & \\
 \hline\hline
 & & \\
 $L4_6$ & \includegraphics{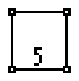} &
 $\displaystyle\frac{[4][6][10]}
 { t^{17} - t^{16} - t^{14} - t^{13} - t^{12} - t^{11} + t^6 + t^5 + t^4 + t^3 + t - 1}$ \cr
 & $(17, 17)$ & \cr
 & & \\
 & & \footnotesize Can be reduced to: \cr
 & & \cr
 & \footnotesize $(11, 11)$ & \footnotesize $\displaystyle\frac{[4][3'][5']}
 {t^{11} - 3t^{10} + 3t^9 - t^8 - 2t^7 + 4t^6 - 4t^5 + 2t^4 + t^3 - 3t^2 + 3t - 1}$ \cr
 & & \\
$L4_7$ & \includegraphics{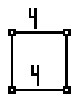} &
 $\displaystyle\frac{[2][4][6]}{t^9 - t^8 - t^7 - t^6 + t^3 + t^2 + t - 1}$ \cr
 & $(9, 9)$  & \cr
 & & \\
 & & \footnotesize Can be reduced to: \cr
 & & \cr
 &  \footnotesize $(7, 7)$  & $\displaystyle\frac{[2][4][3']}{t^7 - 2t^6 + t^4 - t^3 + 2t - 1}$ \cr
 & & \\
 & & \cr
 $L4_8$ & \includegraphics{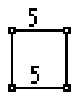} &
 $\displaystyle\frac{[2][6][10]}
 {t^{15} - t^{14} - t^{13} - t^{12} - t^{11} - t^{10} + t^5 + t^4 + t^3 + t^2 + t - 1}$ \cr
 &$(15, 15)$ & \\
 & & \cr
 & & \footnotesize Can be reduced to: \cr
 & & \cr
 & \footnotesize $(9, 9)$  & \footnotesize $\displaystyle\frac{[2][3'][5']}
 { t^9 - 3t^8 + 2t^7 + t^6 - 3t^5 + 3t^4 - t^3 - 2t^2 + 3t - 1}$ \cr
 & & \cr
 $L4_9$ & \includegraphics{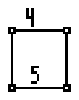} &
 $\displaystyle\frac{[4][6][10]}
 { t^{17} - t^{16} - 2t^{14} - t^{13} - 2t^{12} - t^{10} + t^9 - t^8 + t^7 + 2t^5 + t^4 + 2t^3 + t - 1}$ \cr
 &$(17, 17)$ & \\
 & & \cr
 & & \footnotesize Can be reduced to: \cr
 & & \cr
 & \footnotesize $(11, 11)$ &  $\displaystyle\frac{[4][3'][5']}
    { t^{11} - 3t^{10} + 3t^9 - 2t^8 + 2t^6 - 2t^5 + 2t^3 - 3t^2 + 3t - 1}$ \cr
 & & \\
 \hline\hline
\end{tabular}
 \label{lanner_46_49}
 \end{center}
 \vspace{1mm}
 \end{table*}

\newpage

\begin{table*}
\footnotesize
 \centering
\caption{\hspace{1mm}The Lann\'er diagrams on $5$ vertices and
Poincar\'e series}
 \vspace{1mm}
\begin{tabular}{||c|c|c||}
 \hline\hline
 & & \cr
 $L5_1$ & \includegraphics{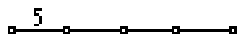} &
 The numerator is $[2][12][20][30]$ \cr
 & & The denominator is \cr
 & $\chi = 1/14400$ & $ t^{60} - t^{59} - t^{53} + t^{52} - t^{51} + t^{50} - t^{49} + t^{46} - t^{45}$ \cr
 & $\deg = (60, 60)$ & $+ t^{44} - 2t^{43} + 2t^{42} - t^{41} + t^{40} - t^{39} + t^{38} - t^{37} + 2t^{36}$ \cr
 & & $- 2t^{35} + 2t^{34} - 2t^{33} + 2t^{32} - t^{31} + t^{30} - t^{29} + 2t^{28} - 2t^{27}$ \cr
 & & $+ 2t^{26} - 2t^{25} + 2t^{24} - t^{23} + t^{22} - t^{21} + t^{20} - t^{19} + 2t^{18}$ \cr
 & & $- 2t^{17} + t^{16} - t^{15} + t^{14} - t^{11} + t^{10} - t^9 + t^8 - t^7 - t + 1$ \cr
\hline\hline
\end{tabular}
 \label{lanner_51}
 \end{table*}

\newpage

\begin{table*}
 \centering
\caption{\hspace{1mm}The quasi-Lann\'er diagrams on $4$ vertices and
Poincar\'e series, none of them reciprocal}
 \vspace{1mm}
 \label{qlanner_ql41_ql49}
\begin{tabular}{||c|c|c|c||}
 \hline\hline
 & & & \cr
 Label & Diagram & Poincar\'e series & Inf.gr.= \cr
 & Degrees& $\chi=0$ in all cases & $b_n + 1$ \cr
 \hline\hline
 & & & \cr
 $QL4_1$ & \includegraphics{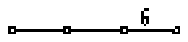} &
$\displaystyle \frac{[2][4][6]}{t^8 - t^4 - t + 1}$ & $1$ \cr
 & $(9,8)$ & &\\
 & & & \cr
 $QL4_2$ & \includegraphics{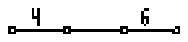} &
 $\displaystyle  \frac{[2][4][6]}{t^8 + t^6 - t^5 - t^3 - t + 1}$ & $1$ \cr
& $(9,8)$ & & \\
 & & & \cr
 $QL4_3$ & \includegraphics{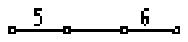} &
 $\displaystyle \frac{[2][6][10]}{t^{14} + t^{12} - t^5 - t^3 - t + 1}$ & $1$ \cr
 & $(15,14)$ & & \\
 & & & \cr
 $QL4_4$ & \includegraphics{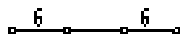} &
 $\displaystyle\frac{[2]^2[6]}{t^7 + t^6 - t^5 - t^2 - t + 1}$ & $2$ \cr
& $(7,7)$ & & \\
 & & & \cr
 $QL4_5$ & \includegraphics{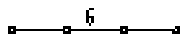} &
 $\displaystyle\frac{[2][6]}{t^6 - t^4 + t^2 - 2t + 1}$ & $2$ \cr
& $(6,6)$ & & \\
 & & & \cr
 $QL4_6$ & \includegraphics{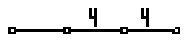} &
 $\displaystyle\frac{[2][4][6]}{t^8 - t^3 - t + 1}$ & $1$ \cr
& $(9,8)$&  &\\
 & & & \cr
 $QL4_7$ & \includegraphics{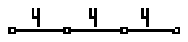} &
 $\displaystyle\frac{[2]^2[4]}{t^5 + t^4 - t^3 - t^2 - t + 1}$ & $2$ \cr
& $(5,5)$ & & \cr
 & & & \cr
 $QL4_8$ & \includegraphics{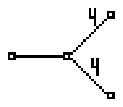} &
 $\displaystyle\frac{[2][4][6]}{t^8 + t^6 - t^5 + t^4 - 2t^3 - t + 1}$ & $1$ \cr
 & $(9,8)$ & & \\
 & & \footnotesize Can be reduced to : & \cr
 &  \footnotesize $(7,6)$ & \footnotesize
 $\displaystyle\frac{[2][4][3']}{t^6 - t^5 + t^4 - t^3 + t^2 - 2t + 1}$  &  \\

 \hline\hline
\end{tabular}
 \end{table*}

\newpage

\begin{table*}
 \centering
\caption{\hspace{1mm}The quasi-Lann\'er diagrams on $4$ vertices and
Poincar\'e series, none of them reciprocal}
 \vspace{1mm}
 \label{qlanner_ql410_ql416}
\begin{tabular}{||c|c|c|c||}
 \hline\hline
 & & & \cr
 Label & Diagram & Poincar\'e series & Inf.gr.= \cr
 & & $\chi=0$ in all cases &  $b_n + 1$ \cr
 \hline\hline
 & & & \\
 $QL4_9$ & \includegraphics{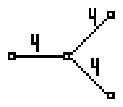} &
 $\displaystyle\frac{[2]^2[4]}{2t^5 + t^4 - 2t^3 - t^2 - t + 1}$ & $3$ \cr
 & $(5,5)$ & & \\
 & & & \cr
 $QL4_{10}$ & \includegraphics{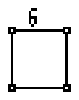} &
 $\displaystyle\frac{[4][6]}{t^8 + t^6 - t^5 - t^3 + t^2 - 2t + 1}$ & $2$ \cr
 & $(8,8)$ & & \\
 & & & \cr
 $QL4_{11}$ & \includegraphics{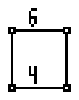} &
 $\displaystyle \frac{[2][4][6]}{t^9 + t^8 + t^7 + t^6 - t^5 - t^4 - t^3 - t^2 - t + 1}$ & $2$ \cr
 &$(9,9)$ & & \\
 & & & \cr
 $QL4_{12}$ & \includegraphics{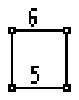} &
$\displaystyle \frac{[2][6][10]} {t^{15} + t^{14} + t^{13} + t^{12}
+ t^{11} - 2t^5 - t^4 - t^3 - t^2 - t + 1}$ & $2$ \cr
 &$(15, 15)$ & & \cr
 & & \footnotesize Can be reduced to: &\\
&&   \footnotesize$\displaystyle \frac{[2][6][5']} {t^{11} - 2t +
1}$ & \cr
 & & & \cr
$QL4_{13}$ & \includegraphics{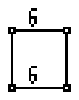} &
 $\displaystyle \frac{[2][6]}{3t^6 - 2t^5 - 2t + 1}$ & $4$ \\
 &$(6,6)$& & \cr
 & & & \cr
 $QL4_{14}$ & \includegraphics{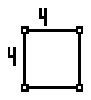} &
 $\displaystyle \frac{[2][4][6]}{t^8 + t^7 - t^3 - t^2 - t + 1}$ & $1$ \\
 &$(9,8)$& & \cr
  & & \footnotesize Can be reduced to: &\\
 &&  \footnotesize$\displaystyle \frac{[2][4][3']}{t^6 - t^4 + t^3 - 2t + 1}$ & \cr
 & & & \cr
 & &  & \cr
 \hline\hline
\end{tabular}
 \end{table*}

\newpage


\begin{table*}
\footnotesize
 \centering
\caption{\hspace{1mm}The quasi-Lann\'er diagrams on $4$ vertices and
Poincar\'e series, none of them reciprocal (cont.)}
 \vspace{1mm}
 \label{qlanner_ql417_ql423}
\begin{tabular}{||c|c|c|c||}
 \hline\hline
 & & & \cr
 Label & Diagram & Poincar\'e series & Inf.gr.= \cr
 & & $\chi=0$ in all cases & $b_n + 1$ \cr
 \hline\hline
 & & & \cr
 $QL4_{15}$ & \includegraphics{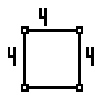} &
 $\displaystyle \frac{[2][4][6]}{(t^9 + t^8 + t^7 - 2t^3 - t^2 - t + 1}$ & $2$ \cr
 &$(9, 9)$ & & \\
 & & \footnotesize Can be reduced to: &\\
 &&  \footnotesize$\displaystyle \frac{[2][4][3']}{t^7 - 2t + 1}$ & \cr
 & & & \cr
 $QL4_{16}$ & \includegraphics{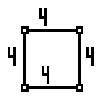} &
 $\displaystyle \frac{[2][4]}{3t^4 - 2t^3 - 2t + 1}$ & $4$ \\
 & $(4,4)$ & & \cr
 $QL4_{17}$ & \includegraphics{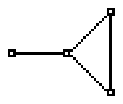} &
 $\displaystyle \frac{[3][4]}{t^4 - t^3 + t^2 - 2t + 1}$ & $1$ \cr
 & $(5,4)$&  & \\
 & &  & \cr
 $QL4_{18}$ & \includegraphics{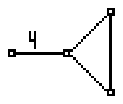} &
 $\displaystyle -\frac{[2][4][6]}{t^8 + t^6 - t^4 - t^2 - t + 1}$ & $1$ \cr
 &$(9, 8)$& & \\
 & & & \cr
 $QL4_{19}$ & \includegraphics{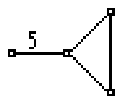} &
 $\displaystyle\frac{[2][6][10]}{t^{14} + t^{12} + t^{11} - t^5 - t^4 - t^2 - t + 1}$ & $1$\cr
 &$(15,14)$ & & \cr
 & & \footnotesize Can be reduced to: &\\
  && \footnotesize$\displaystyle \frac{[2][6][5']}{t^{10} - t^9 + t^8 - t^6 + t^5 - t^4 + t^3 - 2t + 1}$ &  \\
 & & & \cr
 $QL4_{20}$ & \includegraphics{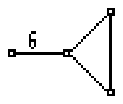} &
 $\displaystyle \frac{[2][6]}{2t^6 - t^5 - t^4 + t^3 - 2t + 1}$ & $3$ \cr
 &$(6,6)$ & & \\
 & & & \cr
 $QL4_{21}$ & \includegraphics{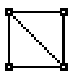} &
 $\displaystyle \frac{[3][4]}{t^5 - 2t + 1}$ & $2$ \cr
 &$(5, 5)$ & &  \cr
 & & & \cr
 $QL4_{22}$ & \includegraphics{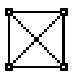} &
 $\displaystyle \frac{[2][3]}{3t^3 - 2t^2 - 2t + 1}$ & $4$ \cr
 &$(3, 3)$ & & \cr
 & & & \\ \hline\hline
\end{tabular}
 \end{table*}

\newpage

\begin{table}
 \centering
\caption{\hspace{1mm}The quasi-Lann\'er diagrams on $5$ vertices and
Poincar\'e series, none of them reciprocal}
 \vspace{1mm}
 \label{qlanner_ql51_ql55}
\begin{tabular}{||c|c|c|c||}
 \hline\hline
 & & & \cr
 Label & Diagram & Poincar\'e series & Inf.gr.= \cr
 & $\chi$, Degrees &  & $b_n + 1$ \\
 \hline\hline
 & & & \cr
 $QL5_1$ & \includegraphics{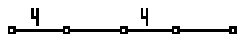} &
 The numerator is $[2][6][8][12]$ & \cr
 & $\chi = -1/1152$& The denominator is & \cr
 & $(24, 23)$ & $t^{23}+t^{19} - t^{18} - t^{16}+t^{15} - 2t^{14}$ & $1$ \cr
 & & $+ t^{13} - 2t^{12}+t^{11} - 2t^{10}+2t^9 - 2t^8 $ & \cr
 & & $+ t^7 - t^6+2t^5 - t^4+t^3+t - 1$ & \cr
  & & & \\

 $QL5_2$ & \includegraphics{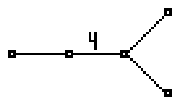} &
 The numerator is $[2][6][8][12]$ & \cr
 & $\chi = -1/576 $ & The denominator is & \cr
 & $(24, 23)$ & $t^{23}+t^{20} - t^{17} - t^{15} - t^{13} - t^{12}$ & $1$ \cr
 & & $- t^{11} - t^8+t^5+t^3+t - 1$ & \cr
  & & & \\

 $QL5_3$ & \includegraphics{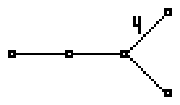} &
 The numerator is $[2][4][5][6][8]$ & \cr
 & $\chi = -1/1920 $ & The denominator is & \cr
 & $(20, 19)$ & $t^{19}+t^{18}+t^{17}+t^{16} - t^{14} -
 t^{13} - 2t^{12} - 2t^{11}$ & $1$ \cr
 & & $- 2t^{10} - t^9 - t^8+t^7+t^6+2t^5+t^4+t^3 - 1$ & \cr
 & & & \\

 $QL5_4$ & \includegraphics{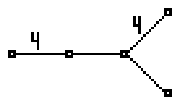} &
 The numerator is $[2][4][6][8]$ & \cr
 & $\chi = -1/384 $& The denominator is & \cr
 & $(16, 16)$ & $ t^{16}+t^{15} - t^{12} - 2t^{10}+t^9 - 3t^8$ & $2$ \cr
 & & $+ t^7 - 2t^6+2t^5 - t^4+2t^3+t - 1$ & \cr
 & & & \\

 $QL5_5$ & \includegraphics{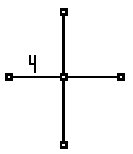} &
 The numerator is $[2][4]^2[6]$ & \cr
 & $\chi = -1/192 $ & The denominator is & \cr
 & $(12, 12)$ & $2t^{12}+t^{11} - 4t^8 - t^7 - 3t^6+2t^5 - t^4+3t^3+t - 1$ & $3$ \cr
 & & & \\
 \hline\hline
\end{tabular}
 \end{table}

\newpage

\begin{table}
 \centering
\caption{\hspace{1mm}The quasi-Lann\'er diagrams on $5$ vertices and
Poincar\'e series, none of them reciprocal (cont.)}
 \vspace{1mm}
 \label{qlanner_ql56_ql59}
\begin{tabular}{||c|c|c|c||}
 \hline\hline
 & & & \cr
 Label & Diagram & Poincar\'e series & Inf.gr.= \cr
 & $\chi$, Degrees &  & $b_n + 1$ \cr
 & & & \cr
 \hline\hline
 & & & \\
 $QL5_6$ & \includegraphics{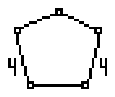} &
 The numerator is $[2][6][8][12]$ & \cr
 & $\chi = -1/144$& The denominator is & \cr
 & $(24, 23)$ & $t^{23}+t^{22}+t^{20} - t^{17} - 2t^{16} - t^{15} -
 t^{14} - t^{13} $ & $1$ \cr
 & & $- 3t^{12} - t^{11} - t^{10} - 3t^8+t^5+t^3+t^2+t - 1$ & \cr
 & & & \\
 $QL5_7$ & \includegraphics{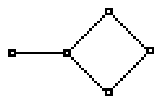} &
 The numerator is $[2][4]^2[5][6]$ & \cr
 & $\chi = -1/960$ & The denominator is & \cr
 & $(16, 15)$ & $t^{15}+t^{14}+t^{13}+t^{12} - t^{11} - 2t^{10} $ & $1$ \cr
 & & $- 2t^9 - 3t^8 - t^7+t^5+2t^4+t^3+t^2 - 1$ & \cr
 & & & \\
 $QL5_8$ & \includegraphics{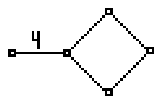} &
 The numerator is $[2][4][6][8]$ & \cr
 & $\chi = -1/192$ & The denominator is & \cr
 & $(16, 16)$ & $t^{16}+t^{15} - t^{11} - t^{10} - t^9 $ & $2$ \cr
 & & $- 2t^8 - t^7+t^3+t^2+t - 1$ & \cr
 & & & \\
 $QL5_9$ & \includegraphics{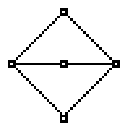} &
 The numerator is $[2][4]^2[6]$ & \cr
 & $\chi = -1/192$ & The denominator is & \cr
 & $(12, 12)$ & $2t^{12}+t^{11} - t^{10}+t^9 - 3t^8 - 2t^7$ & $3$ \cr
 & & $- t^6 - 2t^5+t^3+2t^2+t - 1$ & \cr
 & & & \\
 \hline\hline
\end{tabular}
 \end{table}

\newpage


\begin{table}[h]
 \centering
\caption{\hspace{1mm}The quasi-Lann\'er diagrams on $6$ vertices and
Poincar\'e series}
\begin{tabular}{||c|c|c|c||}
 \hline\hline
 Label & Diagram & Poincar\'e series & Inf.gr.= \cr
 & Degrees & $\chi = 0$ in all cases & $b_n + 1$ \cr
 \hline\hline
 & & & \cr
 $QL6_1$ & \includegraphics{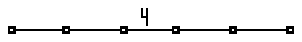} & The numerator is $[2]^2[6][8][12]$ & \cr
 & $(25, 25)$ & The denominator is & \cr
 & & $ t^{25}+t^{24} - t^{23} - t^{22} - t^{19} - t^{18}+2t^{16}$ & $2$ \cr
 & & $- t^{15} - t^{14}+t^{13}+t^{12} - t^{11}+t^{10}+t^7$ & \cr
 & & $+ t^6 - 2t^5+t^3 - t^2 - t+1$ & \cr
\hline\hline
\end{tabular}
 \label{q_lanner_61}
 \end{table}

 \begin{table}[ht]
\footnotesize
 \centering
\caption{\hspace{1mm}The quasi-Lann\'er diagrams on $6$ vertices and
Poincar\'e series (cont.)}
\begin{tabular}{||c|c|c|c||}
 \hline\hline
 Label & Diagram & Poincar\'e series & Inf.gr.= \cr
 & Degrees & $\chi = 0$ in all cases &  $b_n + 1$ \cr
 \hline\hline
 & & & \\
 $QL6_7$ & \includegraphics{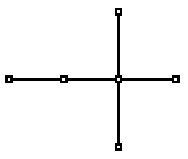} & The numerator is $[2]^2[4][6][8]$ & \cr
 & $(17, 16)$ & The denominator is & \cr
 & & $t^{16} - t^{14} - t^{12}+t^8 - t^7$ & $1$ \cr
 & & $+ 2t^5+t^4 - t^3 - t^2 - t+1$ & \cr
 & & & \\
 $QL6_8$ & \includegraphics{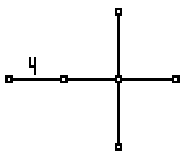} & The numerator is
 $[2]^2[4][6][8]$ & \cr
 & $(17, 17)$ & The denominator is & \cr
 & & $3t^{17}+t^{16} - 3t^{15} - 3t^{14} - t^{13}+t^{12} - 2t^{11}$ & \cr
 & & $+ 4t^8 - t^7+2t^6+2t^4 - 2t^3 - t^2 - t+1$ & $4$ \cr
 & & & \\
 $QL6_9$ & \includegraphics{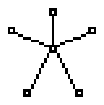} &
 The numerator is $[2]^2[4]^2[6]$ & \cr
 & $(13, 13)$ & The denominator is & \cr
 & & $4t^{13} + t^{12} - 4t^{11} - 6t^{10} - 3t^9 + 3t^8$ & \cr
 & & $+ 5t^6 + 2t^5 + 3t^4 - 4t^3 - t^2 - t + 1$ & \cr
 & & & \\
 & \footnotesize $(11, 11)$ & \footnotesize Can be reduced to numerator: $[2]^2[4]^2[3']$ & \cr
 & & \footnotesize $4t^{11} - 3t^{10} - 5t^9 + 2t^8 + t^6$ & $5$ \cr
 & & \footnotesize$- t^5 + 5t^4 - 2t^3 - 2t + 1$ & \cr
 & & & \\
 \hline\hline
\end{tabular}
 \label{q_lanner_67_69}
 \end{table}

\begin{table}[h]
 \centering
\caption{\hspace{1mm}The quasi-Lann\'er diagrams on $6$ vertices and
Poincar\'e series (cont.)}
 \vspace{1mm}
\begin{tabular}{||c|c|c|c||}
 \hline\hline
 & & & \cr
 Label & Diagram & Poincar\'e series & Inf.gr.= \cr
& $\chi$, Degrees &   & $b_n + 1$ \cr
 & &  & \\
 \hline\hline
 & & & \cr
 $QL6_{10}$ & \includegraphics{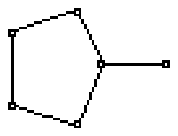} & The numerator is
 $[2][5][6][8]$ & \cr
 & $(17, 16)$ & The denominator is & \cr
 & & $t^{16} - t^{15} - t^{11}+t^9 - t^8+t^7+t^3 - 2t+1$ & $1$ \cr
 & & & \\
 $QL6_{11}$ & \includegraphics{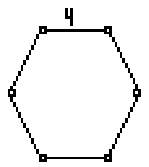} &
 The numerator is $[2][6][8][10][12]$ & \cr
 & $(33, 33)$ & The denominator is & \cr
 & & $t^{33} + t^{32} + t^{31} - t^{28} - 2t^{26} - t^{25} $ & \cr
 & & $- t^{24} -3t^{23} - 3t^{22} - t^{21} - t^{20} - 2t^{19}$ & \cr
 & & $+ t^{18} + 2t^{16} + t^{15} + 2t^{14} + t^{13} + 3t^{12}$ & \cr
 & & $+ t^{11} + t^{10} + t^9 + t^8 - t^5 - t^2 - t + 1$ & \cr
 & & & \\
 & \footnotesize$(29, 29)$ & \footnotesize Can be reduced to numerator:
 $[2][6][8][12][5']$ & \cr
 & & $t^{29} - t^{26}+t^{23} - 2t^{22} - 2t^{19}+t^{18} - t^{15}$ & \cr
 & & $+ t^{14}+2t^{12} - t^{11}+2t^8 - t^6 + t^3 - 2t + 1$ & $2$ \cr
 & & & \\
 $QL6_{12}$ & \includegraphics{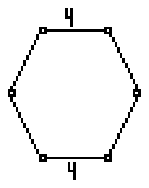} &
 The numerator is $[2][6][8][12]$ & \cr
 & $(24, 24)$ & The denominator is & \cr
 & & $5t^{24} - 4t^{23} - 2t^{21}+3t^{20} - 6t^{19}+3t^{18}$ & \cr
 & & $-6t^{17}+7t^{16} - 4t^{15}+2t^{14} - 6t^{13}+8t^{12} - 4t^{11}$ & \cr
 & & $+ 2t^{10} - 2t^9+5t^8 - 2t^7+3t^6 - 2t^5+t^4 - 2t+1$ & $6$ \cr
 & & & \\
 \hline\hline
\end{tabular}
 \label{q_lanner_610_612}
 \end{table}

\newpage


\begin{table}[h]
 \centering
\caption{\hspace{1mm}The quasi-Lann\'er diagrams on $7$ vertices and
Poincar\'e series}
 \vspace{1mm}
\begin{tabular}{||c|c|c|c||}
 \hline\hline
 & & & \cr
 Label & Diagram & Poincar\'e series & Inf.gr.= \cr
& $\chi$, Degrees &   & $b_n + 1$ \cr
 & &  & \\
 \hline\hline
 $QL7_1$ & \includegraphics{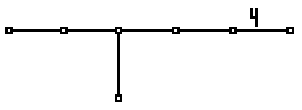} & The numerator is $[2][4][6][8][9][10][12]$ & \cr
 & $(44, 43)$ & The denominator is &  \\
 & &$t^{43} + t^{42} + t^{41} + t^{40} - t^{37} - t^{36} - 2t^{35} - 3t^{34}$ & \cr
 & &$- 3t^{33} - 3t^{32} - 3t^{31} - 2t^{30} - t^{29} + t^{28} + 2t^{27} + 3t^{26}$ & $1$ \cr
 & &$+ 4t^{25} + 4t^{24} + 5t^{23} + 5t^{22} + 4t^{21} + 2t^{20} + t^{19} - t^{18}$ & \cr
 & &$- t^{17} - 3t^{16} - 3t^{15} - 4t^{14} - 3t^{13} - 3t^{12} - t^{11} - t^{10}$ & \cr
 & &$+ t^9 + t^7 + t^6 + t^5 + t^4 + t^3 - 1$ &  \\
 $QL7_2$ & \includegraphics{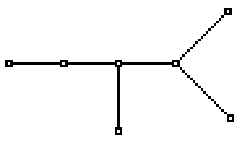} & The numerator is $[2]^2[6][8][9][10][12]$ & \cr
 & $(42,41)$ & The denominator is & \\
 & &$t^{41} + t^{40} - t^{35} - t^{34} - t^{33} - 2t^{32} - 2t^{31} - 2t^{30}$ & \cr
 & &$- t^{29} - t^{28} + 2t^{26} + 2t^{25} + t^{24} + 2t^{23} + 3t^{22} + 4t^{21}$ & $1$ \cr
 & &$+ 3t^{20} + t^{19} - t^{15} - 2t^{14} - 2t^{13} - 2t^{12} - 2t^{11} - t^{10}$ & \cr
 & &$+ t^4 + t^3 + t^2 - 1$ & \cr
 & & & \\
 $QL7_3$ & \includegraphics{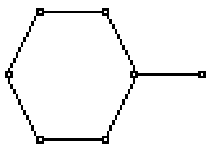} & The numerator is $[2]^2[6][7][8][9][10][12]$ & \cr
 & $(48,47)$ & The denominator is & \\
 & &$t^{47} + 2t^{46} + 2t^{45} + 2t^{44} + 2t^{43} + 2t^{42} + t^{41}$ & \cr
 & &$ - 2t^{40} - 4t^{39} - 6t^{38} - 9t^{37} - 11t^{36} - 12t^{35} - 13t^{34}$ & \cr
 & &$- 11t^{33} - 7t^{32} - 3t^{31} + t^{30} + 6t^{29} + 12t^{28} + 18t^{27}$ & $1$ \cr
 & &$+ 21t^{26} + 22t^{25} + 21t^{24} + 19t^{23} + 17t^{22} + 13t^{21} + 6t^{20}$ & \cr
 & &$- 5t^{18} - 8t^{17} - 10t^{16} - 12t^{15} - 12t^{14} - 10t^{13} - 9t^{12}$ & \cr
 & &$- 6t^{11} - 2t^{10} + t^8 + 2t^7 + 2t^6 + 2t^5 + 2t^4$ & \cr
 & &$+ 2t^3 + t^2 - t - 1$ & \\
 \hline\end{tabular}
 \label{q_lanner_7}
 \end{table}

\newpage


\begin{table}[h]
 \centering
\caption{\hspace{1mm}The quasi-Lann\'er diagrams on $8$ vertices and
Poincar\'e series}
 \vspace{1mm}
\begin{tabular}{||c|c|c||}
 \hline\hline
 & &   \cr
 Label & Diagram & Poincar\'e series   \cr
& Degrees & $\chi = 0$ in all cases  \cr
 & &  \\
 \hline\hline
 $QL8_1$ & \includegraphics{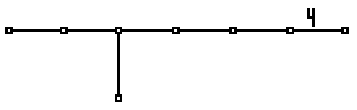} & The numerator is  $[4][6][8][10][12][14][18]$  \cr
 & $(65, 64)$ & The denominator is  \\
 & &$t^{64} + t^{62} - t^{58} - 2t^{56} - 3t^{54} + t^{53} - 3t^{52} + t^{51}$  \cr
 & &$- 3t^{50} + 2t^{49} - 2t^{48} + 3t^{47} - 2t^{46} + 4t^{45} - t^{44} + 4t^{43}$  \cr
\hline\hline
\end{tabular}
 \label{q_lanner_81_82}
 \end{table}

\newpage

\begin{table}[h]
 \centering
\caption{\hspace{1mm}The quasi-Lann\'er diagrams on $8$ vertices and
Poincar\'e series}
 \vspace{1mm}
\begin{tabular}{||c|c|c||}
 \hline\hline
 & &  \cr
 Label & Diagram & Poincar\'e series \cr
& Degrees & $\chi = 0$ in all cases  \cr
 & & \\
 \hline\hline
 $QL8_3$ & \includegraphics{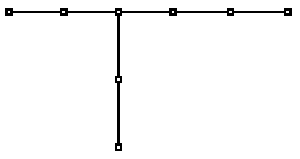} & The numerator is $[2][6][8][10][12][14][18]$  \cr
 & $(63, 62)$ & The denominator is  \\
 & & $t^{62} - t^{56} - t^{54} - t^{52} - t^{50} - t^{49} + t^{48}$ \cr
 & & $- t^{47} + t^{46} + t^{44} + 2t^{42} + t^{40} + 2t^{39} + t^{37}$ \cr
\hline\hline
\end{tabular}
 \label{q_lanner_83_84}
 \end{table}

\newpage


\begin{table}[h]
 \centering
\caption{\hspace{1mm}The quasi-Lann\'er diagrams on $9$ vertices and
Poincar\'e series}
 \vspace{1mm}
\begin{tabular}{||c|c|c||}
 \hline\hline
 & &  \cr
 Label & Diagram & Poincar\'e series  \cr
& $\chi$, Degrees &   \cr
 & &  \\
 \hline\hline
 & &  \cr
 $QL9_1$ & \includegraphics{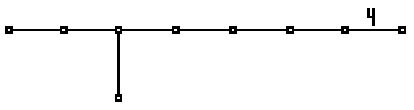} & The numerator is [2][12][14][16][18][20][24][30]   \cr
 & $(128, 127)$ & The denominator is \\
 & & $t^{127} - t^{117} + t^{116} - 2t^{115} + t^{114} - 3t^{113}$ \cr
\hline\hline
\end{tabular}
 \label{q_lanner_91}
 \end{table}

\newpage


\begin{table*}[h]
 \centering
\caption{\hspace{1mm}The quasi-Lann\'er diagrams on $9$ vertices and
Poincar\'e series (cont.)}
 \vspace{1mm}
\begin{tabular}{||c|c|c||}
 \hline\hline
 & &  \cr
 Label & Diagram         & Poincar\'e series  \cr
       & $\chi$, Degrees &  \cr
 & &  \\
 \hline\hline
  & &  \cr
 $QL9_3$ & \includegraphics{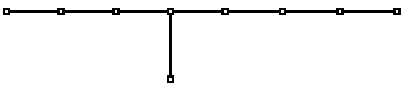} &
    The numerator is $[2][8][12][14][18][20][24][30]$  \cr
 & $(120, 119)$ & The denominator is   \\
 & &  $t^{119} - t^{111} - t^{107} - t^{105} - t^{103} + t^{102} - 2t^{101}$   \cr
\hline\hline
\end{tabular}
 \label{q_lanner_9}
 \end{table*}

\newpage


\begin{table*}[h]
 \footnotesize
 \centering
\caption{\hspace{1mm}The quasi-Lann\'er diagrams on $10$ vertices
and Poincar\'e series}
 \vspace{1mm}
\begin{tabular}{||c|c|c||}
 \hline\hline
 & & \cr
 Label & Diagram & Poincar\'e series \cr
& Degrees & $\chi = 0$ in all cases \cr
& &  \cr
 \hline\hline
& &  \cr
 $QL10_1$ & \includegraphics{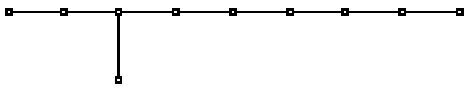} &
 The numerator is $[2]^2[12][14][16][18][20][24][30]$\cr
 & $(129, 128)$ & The denominator is \cr
 & & $t^{128} - t^{126} - t^{116} + t^{105} - t^{103} + t^{102} + t^{101}$ \cr
\hline\hline
\end{tabular}
 \label{q_lanner_9_10}
 \end{table*}

\clearpage

\begin{table*}[h] 
 \centering
\caption{\hspace{1mm}The quasi-Lann\'er diagrams on $10$ vertices
and Poincar\'e series (cont.)}
 \vspace{1mm}
\begin{tabular}{||c|c|c||}
 \hline\hline
 & & \cr
 Label & Diagram & Poincar\'e series \cr
& Degrees & $\chi = 0$ in all cases \cr
 \hline\hline$QL10_3$ & \includegraphics{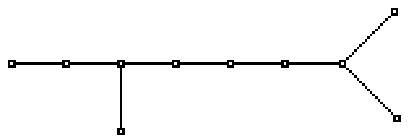} & The numerator is $[2]^2[12][14][16][18][20][24][30]$ \cr
 & $(129, 129)$ & The denominator is \cr
 & &$2t^{129} + t^{128} - 2t^{127} - t^{126} - 2t^{122} + 2t^{121} + 2t^{120}$ \cr
 & &$- 2t^{119} - 3t^{118} + t^{117} - t^{116} - 3t^{114} + 3t^{113} + 2t^{111}$ \cr
 & &$- 6t^{110} + 5t^{109} - t^{108} + 2t^{107} - 7t^{106} + 12t^{105} - 3t^{104}$ \cr
 & &$+ 5t^{103} - 10t^{102} + 13t^{101} - 6t^{100} + 10t^{99} - 16t^{98} + 18t^{97}$ \cr
 & The number of infinite &$- 5t^{96} + 10t^{95} - 21t^{94} + 23t^{93} - 14t^{92} + 12t^{91} - 22t^{90}$ \cr
 & special subgroups: $3$ &$+ 24t^{89} - 14t^{88} + 14t^{87} - 29t^{86} + 23t^{85} - 13t^{84} + 10t^{83}$ \cr
 & &$- 31t^{82} + 26t^{81} - 13t^{80} + 8t^{79} - 24t^{78} + 16t^{77} - 10t^{76}$ \cr
 & &$+ 5t^{75} - 22t^{74} + 14t^{73} + t^{72} - 5t^{71} - 13t^{70} + 9t^{69}$ \cr
 & &$+ 7t^{68} - 13t^{67} - t^{65} + 19t^{64} - 17t^{63} + 6t^{62} - 9t^{61}$ \cr
 & &$+ 30t^{60} - 27t^{59} + 16t^{58} - 13t^{57} + 36t^{56} - 33t^{55} + 26t^{54}$ \cr
 & &$- 21t^{53} + 37t^{52} - 34t^{51} + 26t^{50} - 24t^{49} + 44t^{48} - 38t^{47}$ \cr
 & &$+ 26t^{46} - 24t^{45} + 38t^{44} - 37t^{43} + 25t^{42} - 24t^{41} + 34t^{40}$ \cr
 & &$- 31t^{39} + 21t^{38} - 24t^{37} + 29t^{36} - 28t^{35} + 16t^{34} - 17t^{33}$ \cr
 & &$+ 22t^{32} - 22t^{31} + 13t^{30} - 14t^{29} + 16t^{28} - 17t^{27} + 9t^{26}$ \cr
 & &$- 9t^{25} + 13t^{24} - 11t^{23} + 5t^{22} - 5t^{21} + 8t^{20} - 8t^{19}$ \cr
 & &$+ 4t^{18} - 3t^{17} + 6t^{16} - 4t^{15} + 2t^{14} - t^{13} + 3t^{12}$ \cr
 & &$- 3t^{11} + 3t^8 - 2t^7 + t^4 - t^2 - t + 1$ \\
 \hline\hline
\end{tabular}
 \label{q_lanner_9_10a}
 \end{table*} 

\clearpage

\clearpage

\newpage

\clearpage

\begin{table*}[h]
 \begin{center}
\caption{\hspace{1mm}The multiparameter Poincar\'e series of affine
Coxeter groups}
 \vspace{1mm}
 \renewcommand{\arraystretch}{1.7}
\begin{tabular}{||c|c|c||}
 \hline\hline
Affine Coxeter & Extended Dynkin  & Poincar\'e  \cr
 group & diagram & series \\
 \hline\hline
 & & \\
  $B_n^{(1)}$ &\includegraphics{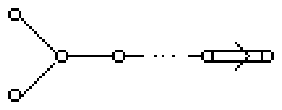} & $\displaystyle\frac{1-t_1^n}{(1-t_1)^n}\left(
  \mathop{\prod}\limits_{i=1}^{n-1}(1-t_1^{2i})\right)
  \mathop{\prod}\limits_{i=0}^{n-1}\frac{1-t_1^{i}t_2}{1-t_1^{n-1+i}t_2}$\\
 $C_n^{(1)}$ &\includegraphics{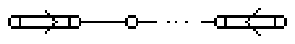} & $\displaystyle\mathop{\prod}\limits_{i=0}^{n-1}
 \frac{(1-t_1^{i+1})(1+t_1^{i}t_2)(1+t_1^{i}t_3)}{(1-t_1)(1-t_1^{n-1+i}t_2t_3)}$\\
 & & \\
 $F_4^{(1)}$ &\includegraphics{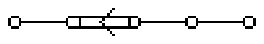} & $\begin{array}{l}\displaystyle\frac{(1+t_1t_2^2)(1+(t_1t_2)^2)
 (1+(t_1t_2)^{3})}
 {(1-t_1^3t_2)(1-t_1^4t_2^3)(1-t_1^5t_2^3)(1-t_1^6t_2^5)}\times\\
 \displaystyle\mathop{\prod}\limits_{i=1}^{3}
 \frac{(1-t_1^{i+1})(1+t_1^{i}t_2)(1-t_2^{i})}{(1-t_1)(1-t_2)}\end{array}$\\
 $G_2^{(1)}$ &\includegraphics{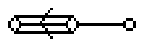} & $\displaystyle\frac{(1+t_1)(1+t_1+t_1^2)(1+t_2)
 (1+t_1t_2+(t_1t_2)^{2})}
 {(1-t_1^2t_2)(1-t_1^3t_2^2)}$ \\
\hline\hline
 \end{tabular}
 \label{tab_multip}
 \end{center}
\end{table*}

\begin{table*}[h]
 \begin{center}
\caption{\hspace{1mm}The real poles and the extremal absolute values
of the non-real poles of the Poincar\'e series. Lann\'er cases}
 \vspace{1mm}
 \renewcommand{\arraystretch}{1.7}
\begin{tabular}{||c|c|c|c|c||}
 \hline\hline
$L5_1$&$L5_2$&$L5_3$&$L5_4$&$L5_5$\\
\hline$\begin{array}{l}
0.833415\\
0.94166\\
1.06195\\
1.19988\\
m = 0.97149\\
M= 1.02935
\end{array}$&$
\begin{array}{l}
0.720106\\
0.898971\\
1.11238\\
1.38868\\
m = 0.96401\\
M= 1.03734\end{array}$&$
\begin{array}{l}
0.659358\\
0.875566\\
1.14212\\
1.51663\\
m = 0.93176\\
M= 1.07324\end{array}$&$
\begin{array}{l}
\text{yes, correct:}\\
\text{no}\\
\text{real}\\
\text{roots}\\
m = 0.94718\\
M = 1.05577\end{array}$&$
\begin{array}{l}
0.61621\\
0.85384\\
1.17118\\
1.62282\\
m = 0.89454\\
M = 1.11788\end{array}$\\
\hline\hline\end{tabular}

 \label{Tabpoles}
 \end{center}
\end{table*}

\begin{table*}[h]
 \begin{center}
\caption{\hspace{1mm}The real poles and the extremal absolute values
of the non-real poles of the Poincar\'e series. Quasi-Lann\'er cases
(The trivial pole 1 is not indicated)}
 \vspace{1mm}
 \renewcommand{\arraystretch}{1.7}\begin{tabular}{||c|c|c|c|c||}
 \hline\hline
$QL4_1$&$QL4_2$&$QL4_3$&$QL4_4$&$QL4_5$\\
\hline$\begin{array}{l}
0.771327\\
m=0.950357\\
M=1.103357
\end{array}$&$
\begin{array}{l}
\\
m = 1\\
M=1.210606
\end{array}$&$
\begin{array}{l}
0.639025
\\
m=0.960217
\\
M=1.142917
\end{array}$&$
\begin{array}{l}
-1.61803
\\
0.618034
\\
m =
M = 1\end{array}$&$
\begin{array}{l}
0.667961
\\
m=0.910638\\
M=1.343628
\end{array}$\\
\hline\hline\end{tabular}

\begin{tabular}{||c|c|c|c|c||}
 \hline\hline
$QL4_6$&$QL4_7$&$QL4_8$&$QL4_9$&$QL4_{10}$\\
\hline$
\begin{array}{l}
0.708134
\\
m=0.957686
\\
M =1.146305
\end{array}$&$\begin{array}{l}
-1.61803
\\
0.618034
\\
m = M= 1
\end{array}$&$
\begin{array}{l}
0.636883
\\
m=1.099895
\\
M=1.139254
\end{array}$&$
\begin{array}{l}
-1.33552
\\
0.552965
\\
m=M=0.822833
\end{array}$&$
\begin{array}{l}
0.561856
\\
m=0.909844\\
M=1.287859\end{array}$\\
\hline\hline\end{tabular}

\begin{tabular}{||c|c|c|c|c|c||}
 \hline\hline
$QL4_{11}$&$QL4_{12}$&$QL4_{13}$&$QL4_{14}$&$QL4_{15}$\\
\hline
$\begin{array}{l}
-1.29065\\
0.51879\\
m=1
\\
M=1.222085
\end{array}$&$
\begin{array}{l}
-1.11231
\\
0.500245
\\
m=1.032895
\\
M=1.107983
\end{array}$&$\begin{array}{l}
\\
0.492432\\
m=0.902209
\\
M=0.911924
\end{array}$&$
\begin{array}{l}
\\
0.551753
\\
m=1.076010\\
M=1.251157
\end{array}$&$
\begin{array}{l}
-1.19004\\
0.504138\\
m=1.103491\\
M=1.169974
\end{array}$\\
\hline\hline\end{tabular}
 \end{center}
\end{table*}

\begin{table*}[h]
 \begin{center}

\begin{tabular}{||c|c|c|c||}
 \hline\hline
$QL4_{16}$&$QL4_{17}$&$QL4_{18}$&$QL4_{19}$\\
\hline$
\begin{array}{l}
\\
0.469396
\\
m=M=0.842693
\end{array}$&$
\begin{array}{l}
\\
0.682328
\\
m=M=1.210606
\end{array}$&$
\begin{array}{l}
0.588985
\\
m=0.962999\\
M=1.231827
\end{array}$&$\begin{array}{l}
0.552531
\\
m=0.986410\\
M=1.243136
\end{array}$\\
\hline\hline\end{tabular}
\begin{tabular}{||c|c|c||}
 \hline\hline
$QL4_{20}$&$QL4_{21}$&$QL4_{22}$\\
\hline$
\begin{array}{l}
0.537613
\\
m=0.942397\\
M=1.023332
\end{array}$&$
\begin{array}{l}
\\
-1.29065\\
0.51879
\\
m=M=1.222085
\end{array}$&$
\begin{array}{l}
-0.767592\\
0.434259
\\
m = 0.94718\\
M = 1.05577\end{array}$\\
\hline\hline\end{tabular}

 \label{Tabpoles4}
 \end{center}
\end{table*}

\begin{table*}[h]
 \begin{center}
\caption{\hspace{1mm}The real poles and the extremal absolute values
of the non-real poles of the Poincar\'e series. Quasi-Lann\'er
cases}
 \vspace{1mm}
 \renewcommand{\arraystretch}{1.7}
\begin{tabular}{||c|c|c|c|c||}
 \hline\hline
$QL5_1$&$QL5_2$&$QL5_3$&$QL5_4$&$QL5_5$\\
\hline$\begin{array}{l}
\\
0.698956\\
0.891273\\
1.09813
\\
m=0.957885\\
M=1.154830
\end{array}$&$
\begin{array}{l}
-1.236\\
-1.05414\\
0.654741\\
0.872516\\
1.12047\\
m=0.914538\\
M=1.161598\end{array}$&$
\begin{array}{l}
\\
0.72899\\
0.903396\\
1.08431
\\
m=0.948289\\
M=1.150045
\end{array}$
&$
\begin{array}{l}
-1.62934\\
0.627864\\
0.862852\\
1.12033
\\
m=0.915917\\
M=1.131353
\end{array}$
&$
\begin{array}{l}
-1.46751\\
0.579431\\
0.842435\\
1.1368
\\
m=0.849730\\
M=1.074701
\end{array}$\\
\hline\hline\end{tabular}

\begin{tabular}{||c|c|c|c||}
 \hline\hline
$QL5_6$&$QL5_7$&$QL5_8$&$QL5_9$\\
\hline$
\begin{array}{l}
\\
0.537456\\
0.820844\\
1.18562
\\
m=0.903486\\
M=1.344774
\end{array}$&$
\begin{array}{l}
\\
0.662566\\
0.876238\\
1.11333
\\
m=0.947039\\
M=1.119572\
\end{array}$
&$
\begin{array}{l}
\\
-1.16298
0.55887\\
0.831791\\
1.15607
\\
m=0.920577\\
M=1.106633
\end{array}$&$
\begin{array}{l}
-1.38639\\
0.491695\\
0.800368\\
1.18595\\
m=0.864826\\
M=1.082570
\end{array}$\\
\hline\hline\end{tabular}

 \label{Tabpoles5}
 \end{center}
\end{table*}

\begin{table*}[h]
 \begin{center}
\caption{\hspace{1mm}The real poles and the extremal absolute values
of the non-real poles of the Poincar\'e series. Quasi-Lann\'er
cases. (The trivial pole 1 is not indicated)}
 \vspace{1mm}
 \renewcommand{\arraystretch}{1.7}\begin{tabular}{||c|c|c|c|c||}
 \hline\hline
$QL6_1$&$QL6_2$&$QL6_3$&$QL6_4$&$QL6_5$\\
\hline$\begin{array}{l}
-1.41222\\
0.741226\\
0.864041\\
\\
m=0.947515\\
M=1.137513
\end{array}$&$
\begin{array}{l}
0.801198\\
0.896819\\
1.07859\\
\\
m=0.966866\\
M=1.130358
\end{array}$&$
\begin{array}{l}
\\
\\
\\
\\
m=0.9459183\\
M=1.136787
\end{array}$&$
\begin{array}{l}
-1.35548\\
0.667522\\
0.822211\\
1.13349
\\
m=0.940570\\
M=1.068224
\end{array}$&$
\begin{array}{l}
0.744209\\
0.8651\\
1.10943\\
\\
m=0.961618\\
M=1.147472
\end{array}$\\
\hline\hline\end{tabular}

\begin{tabular}{||c|c|c|c|c||}
 \hline\hline
$QL6_6$&$QL6_7$&$QL6_8$&$QL6_9$&$QL6_{10}$\\
\hline$
\begin{array}{l}
-1.30069\\
0.634641\\
0.802625\\
1.15397
\\
m=0.915468\\
M=1.051409
\end{array}$&$\begin{array}{l}
0.702245\\
0.840655\\
1.13578\\
\\
m=0.881099\\
M=1.177095
\end{array}$&$
\begin{array}{l}
-1.28976\\
0.59287\\
0.777193\\
1.17396
\\
m=0.897402\\
M=0.981604
\end{array}$&$
\begin{array}{l}
-1.25535\\
0.542596\\
0.744994\\
1.20423
\\
m=0.747308\\
M=0.932831
\end{array}$&$
\begin{array}{l}
0.657119\\
0.814442\\
1.1665\\
\\
m=0.993132\\
M=1.106255
\end{array}$\\
\hline\hline\end{tabular}

\begin{tabular}{||c|c||}
 \hline\hline
$QL6_{11}$&$QL6_{12}$\\
\hline$
\begin{array}{l}
-1.06958\\
0.604368\\
0.784304\\
1.17725
\\
m=0.907820\\
M=1.237682
\end{array}$&$
\begin{array}{l}
0.533802\\
0.742603\\
1.19555\\
\\
m=0.889811\\
M=1.010751
\end{array}$\\
\hline\hline\end{tabular}

 \label{Tabpoles6}
 \end{center}
\end{table*}

\begin{table*}[h]
 \begin{center}
\caption{\hspace{1mm}The real poles and the extremal absolute values
of the non-real poles of the Poincar\'e series. Quasi-Lann\'er
cases}
 \vspace{1mm}
 \renewcommand{\arraystretch}{1.7}
\begin{tabular}{||c|c|c||}
 \hline\hline
$QL7_1$&$QL7_2$&$QL7_3$\\
\hline$\begin{array}{l}
0.755431\\
0.847161\\
0.94682\\
1.0582\\
1.11047
\\
m=0.943848\\
M=1.145506
\end{array}$&$
\begin{array}{l}
0.726032\\
0.827211\\
0.939344\\
1.06657\\
1.13199
\\
m=0.937174\\
M=1.141785
\end{array}$&$
\begin{array}{l}
0.656509\\
0.780576\\
0.92162\\
1.08703\\
1.18705
\\
m=0.892428\\
M=1.160250
\end{array}$\\
\hline\hline\end{tabular}

 \label{Tabpoles7}
 \end{center}
\end{table*}

\begin{table*}[h]
 \begin{center}
\caption{\hspace{1mm}The real poles and the extremal absolute values
of the non-real poles of the Poincar\'e series. Quasi-Lann\'er
cases. (The trivial pole 1 is not indicated)}
 \vspace{1mm}
 \renewcommand{\arraystretch}{1.7}\begin{tabular}{||c|c|c|c||}
 \hline\hline
$QL8_1$&$QL8_2$&$QL8_3$&$QL8_4$\\
\hline$\begin{array}{l}
-1.25799\\
-1.04243\\
0.763804\\
0.837519\\
0.915858
\\
m=0.939234\\
M=1.100237
\end{array}$&$
\begin{array}{l}
-1.29534\\
-1.0366\\
0.741982\\
0.821394\\
0.907002
\\
m=0.921189\\
M=1.114182
\end{array}$&$
\begin{array}{l}
0.77866\\
0.84753\\
0.92114\\
\\
\\
m=0.932965\\
M=1.134982
\end{array}$&$
\begin{array}{l}
0.657583\\
0.760101\\
0.873128\\
\\
\\
m=0.926349\\
M=1.171854
\end{array}$\\
\hline\hline\end{tabular}

 \label{Tabpoles8}
 \end{center}
\end{table*}

\begin{table*}[h]
 \begin{center}
\caption{\hspace{1mm}The real poles and the extremal absolute values
of the non-real poles of the Poincar\'e series. Quasi-Lann\'er
cases}
 \vspace{1mm}
 \renewcommand{\arraystretch}{1.7}\begin{tabular}{||c|c|c|c||}
 \hline\hline
$QL9_1$&$QL9_2$&$QL9_3$&$QL9_4$\\
\hline$\begin{array}{l}
-1.28534\\
-1.02229\\
0.770075\\
0.831719\\
0.896113\\
0.964269\\
1.03702
\\
m=0.935688\\
M=1.108580
\end{array}$&$
\begin{array}{l}
-1.23784\\
-1.0607\\
0.753304\\
0.818656\\
0.88767\\
0.96124\\
1.04029
\\
0.913000\\
1.121820
\end{array}$&$
\begin{array}{l}
-1.19828\\
-1.01676\\
0.826149\\
0.873169\\
0.922226\\
0.973477\\
1.02722
\\
m=0.958697\\
M=1.093674
\end{array}$&$
\begin{array}{l}
-1.03174\\
0.659124\\
0.746793\\
0.841212\\
0.944445\\
1.05875\\
\\
m=0.878261\\
M=1.347490
\end{array}$\\
\hline\hline\end{tabular}

 \label{Tabpoles9}
 \end{center}
\end{table*}

\begin{table*}[h]
 \begin{center}
\caption{\hspace{1mm}The real poles and the extremal absolute values
of the non-real poles of the Poincar\'e series. Quasi-Lann\'er
cases. (The trivial pole 1 is not indicated)}
 \vspace{1mm}
 \renewcommand{\arraystretch}{1.7}
\begin{tabular}{||c|c|c||}
 \hline\hline
$QL10_1$&$QL10_2$&$QL10_3$\\
\hline$\begin{array}{l}
-1.14077\\
-1.008\\
0.878674\\
0.907888\\
0.93783\\
0.968518\\
1.03182
\\
m=0.968106\\
M=1.068399
\end{array}$&$
\begin{array}{l}
-1.01583\\
0.774744\\
0.827863\\
0.882622\\
0.939841\\
1.06238\\
\\
m=0.932746\\
M=1.342964
\end{array}$&$
\begin{array}{l}
-1.014947\\
0.7615867\\
0.8172635\\
0.8751652\\
0.9359320\\
1.0662355\\
\\
m = 0.933795\\
M=1.219335\end{array}$\\
\hline\hline\end{tabular}

 \label{Tabpoles10}
 \end{center}
\end{table*}

\clearpage


\label{lastpage}

\end{document}